\documentclass[12pt]{article}

\usepackage{amsmath, amsthm, amssymb, amsfonts, enumerate}
\usepackage{graphicx}
\usepackage{wasysym, pdfpages}
\usepackage{epstopdf}
\usepackage[colorlinks=true,linkcolor=blue,urlcolor=blue]{hyperref}

\usepackage{geometry}
\def \R{I\!\!R}
\def \E{I\!\!E}
\def \P{I\!\!P}
\def \N{I\!\!N}
\def\two{I\!I}

\def\p{\vskip4truept \noindent}

\newtheorem{thm}{Theorem}[section]

\newtheorem{lem}[thm]{Lemma}
\newtheorem{pro}[thm]{Proposition}
\newtheorem{defi}[thm]{Definition}
\newtheorem{rem}[thm]{Remark}
\newtheorem{nota}[thm]{Notation}

\newtheorem{cla}[thm]{Claim}

\def \conv{ {\rm conv} }

\def \conv{ {\rm conv} }
\def \cav{ {\rm cav} }
\def \vex{ {\rm vex} }
\def \Max{ {\rm Max} }
\def \Min{ {\rm Min} }

\def \lim{ {\rm lim} }
\def \lims{ {\rm limsup} }
\def \limi{ {\rm liminf} }

\let\inf\relax \DeclareMathOperator*\inf{\vphantom{p}inf}

\def\abstract{\begin{center} \small\bf Abstract\end{center}\small}

\title{The value of Markov Chain Games with incomplete information on both sides}
 \author{Fabien Gensbittel\thanks{TSE (GREMAQ, Universit\' e Toulouse 1 Capitole), 21 all\' ee de Brienne, 31000
Toulouse, France. E-mails: \href{mailto:fabien.gensbittel@tse-fr.eu}{fabien.gensbittel@tse-fr.eu}, \href{mailto:jerome.renault@tse-fr.eu}{jerome.renault@tse-fr.eu}} , J\'er\^ome Renault$^*$}
 
\date{\today}
\begin{document}     
\maketitle

\abstract We consider  zero-sum repeated games with incomplete information on both sides, where the states privately observed by each player  follow  independent Markov chains. It  generalizes the model, introduced by Aumann  and  Maschler in the sixties and solved by Mertens and Zamir in the seventies, where the private states of the players were fixed. It also includes  the model introduced in Renault  \cite{R2006}, of Markov chain repeated games with lack of information on  one side, where only one player privately observes the sequence of states. We prove here that the limit value exists, and we obtain  a characterization via the Mertens-Zamir system, where the ``non revealing value function" plugged in the system is now defined as the limit value of an auxiliary ``non revealing" dynamic  game. This non revealing game is defined by restricting the players not to reveal any information on the {\it limit behavior} of their own Markov chain, as in Renault 2006. There are two  key technical difficulties in the proof: 1)  proving regularity, in the sense of equicontinuity, of the $T$-stage non revealing value functions, and 2) constructing strategies by blocks in order to link the values of the non revealing games with the original values.\\

\noindent {\it Key words:} repeated games; incomplete information; zero-sum games; Markov chain; lack of information on both sides;
stochastic games; Mertens-Zamir system.

\section{Introduction}
The pioneering work of Aumann and Maschler in the 1960s introducing  Repeated games with incomplete information (see \cite{aumasch} for a re-edition of their work) was widely studied and extended. One of their most famous results concerns the model of two-player zero-sum repeated game with lack of information on one side and perfect observation: At the beginning of the game, a state variable $k$ is chosen at random in a finite set $K$ using some probability $p$, and announced to Player 1 only. A finite zero-sum game $G_k$ depending on the state variable is then played repeatedly and after each stage the actions played are observed by both players. The famous $\cav u$ theorem of Aumann and Maschler then reads as follows:  the value of the infinitely repeated game exists and is  characterized as  the concave hull of the value of an auxiliary game called the non-revealing game. In a previous work of Renault \cite{R2006}, this model was extended in the following direction: instead of being fixed once for all at the beginning of the game, the state variable was assumed  to evolve according to a Markov chain of initial probability $p$ and transition matrix $M$ over $K$. Renault proved that in this model, called Markov chain game with incomplete information on one side,  the value of the infinitely repeated game exists and gave a characterization for it, very close  to the one of Aumann and Maschler. The main difference was in the definition of the non-revealing game: it was defined in Renault \cite{R2006} as a dynamic game with infinite length where Player 1 was constrained not to reveal any information on the asymptotic behavior of the Markov chain rather than on the current state variable itself. The model of Markov chain game with incomplete information on one side was also studied by Neyman \cite{neyman} who proved the existence of an optimal strategy for Player 1 and also proved a generalization of the existence of the value in case of imperfect observation of actions.    

A more involved model also introduced in \cite{aumasch} is the model of two-player zero-sum repeated game with lack of information on both sides and perfect observation: At the beginning of the game, two states of nature $k$ and $l$ are chosen independently in some finite sets $K$ and $L$, according to some probabilities $p$ and $q$. The state $k$ is announced to Player 1 only, and the state $l$ is announced to Player 2 only.
These states of nature determine a finite zero-sum game $G_{k,l}$ which is then played repeatedly; after each stage the
actions played are observed by both players. Aumann and Maschler showed that the infinitely repeated game may have no value.
Mertens and Zamir \cite{mertenszamir}, showed the existence of the limit value (i.e. the existence of a limit for the sequence of values of finitely repeated games when the number of repetitions goes to infinity) in this model and gave a characterization for it (see the system  \eqref{systemeMZ} later) based on a system of functional equations acting on the value of the non-revealing game.   

The results of Aumann and Maschler and of Mertens and Zamir led to a great number of works dedicated to generalizations of this model, or close extensions of it. Let us cite for games with incomplete information on one side, Kohlberg \cite{kohlberg} for an explicit construction of an
optimal strategy for Player 2, De Meyer and Rosenberg \cite{bdm} and Laraki \cite{ridadual} for alternative proofs of existence of the limit value based on duality, Gensbittel \cite{fg2012} for the extension to infinite action spaces, or  Sorin \cite{sorin83}, Hart \cite{hart},   Simon {\it et al.} \cite{simon}, Renault \cite{R2000} \cite{R2001} for the nonzero-sum case.

Concerning incomplete information on both sides, let us mention the proof of existence of asymptotically optimal strategies by Heuer \cite{heuer},  the   extension of the Mertens-Zamir system by Sorin \cite{sorin84},  the study of an abstract game called the splitting game related to the system of functional equations by Laraki \cite{ridasplitting}, alternative proofs of existence of the limit value based on the so-called operator approach by Rosenberg and Sorin \cite{operator}, or more recently on continuous-time approach by Cardaliaguet, Laraki and Sorin \cite{CLSsiam}.

We only consider here two-player zero-sum games and generalize the model of repeated game with lack of information on both sides to the case where the states of nature $(k,l)$ are no longer fixed at the beginning of the game, but evolve  according to given independent Markov chains $(k_t)_{t\geq 1}$ and $(l_t)_{t \geq 1}$. At the beginning of each stage, $k_t$ is observed by Player 1 only, and $l_t$ is observed by Player 2 only. 
We call such games Markov chain games with lack of information on both sides. 
Note that this model admits as a special case the model of repeated game with incomplete information on both sides, in which the value of the infinitely repeated game may not exist. In this paper, we generalize both the proofs of Renault \cite{R2006} and Mertens and Zamir \cite{mertenszamir} and show the existence of the limit value for Markov chain games with lack of information on both sides. We also give a characterization for this limit value based on a system of functional equations, similar to the one introduced by Mertens and Zamir, and on the generalized notion of non-revealing games as introduced by Renault.

However, as already mentioned in these two works, our expression of the value cannot be easily computed from the basic data of the game. It was already noticed in \cite{mertenszamir}, where several examples of computations are given, underlying that the limit value may be outside of the class of semi-algebraic functions, which contrasts with properties of the limit value for stochastic games  (with complete  information and finite state and action spaces). Moreover, the problem of computation of the value of non-revealing games for Markov chain games with incomplete information on one side was already mentioned by Renault who gave an example which appeared to be very difficult to compute, except for very particular values of the transition matrix   as shown in the work of H\"{o}rner {\it et al.} \cite{horner}.

Section 2 of this paper contains the model. Section 3 contains preliminary results and notations. A few important examples, and the main ideas of the proof, are presented in section 4. In section 5, we define and study a notion of nonrevealing strategies via projection matrices that will be used to quantify a notion of relevant information for both players. We introduce the auxiliary games called nonrevealing games, where both players are restricted to play a nonrevealing strategy, and prove that they have a limit value. 
In section 6, we study some properties of the system of functional equations introduced by Mertens and Zamir associated to the limit value of the non-revealing game. Our main result theorem \ref{thm1} will imply that the limit value is the unique solution of this system.
Section 7 contains a list of open questions related to the possible extensions of the model. Section 8 is an appendix containing the main proofs and some technical results.
 
\section{Model}

Given a finite set $S$, $\Delta(S)$ denotes the set of probabilities over $S$ and the cardinal of $S$ is, with a slight abuse of notations, also denoted by $S$. The set of positive integers is denoted $\N^*$. 

We consider  a zero-sum game between 2 players with action sets $I$ and $J$, sets of states $K$ and $L$, where $I$, $J$, $K$, $L$ are disjoint finite non empty sets, and payoff function 
\[ g:K \times L \times I  \times J \longrightarrow [-1,1] . \]
$M$ and $N$ are given  Markov matrices on $K$ and $L$ respectively, i.e. $M=(M_{k,k'})_{k,k' \in K}$ is a $K\times K$ matrix with non-negative entries and such that for all $k \in K$, $\sum_{k' \in K} M_{k,k'}=1$, and similarly for $N$. The letters $p\in  \Delta(K)$ and $q\in \Delta(L)$ denote initial probabilities.

The game is played by stages in discrete time. We are given two independent Markov chains:  $(k_t)_{t\geq 1}$ with initial distribution $p$ and transition matrix $M$, and $(l_t)_{t\geq 1}$ with initial distribution $q$ and transition matrix $N$. At the beginning of every stage $t\geq 1$,    
Player 1 observes $k_t$ and Player 2 observes $l_t$. Then both players simultaneously  select an action in their action set, if $(i_t,j_t)$ in $I\times J$ is played then Player 1's payoff for stage $t$ is $g(k_t,l_t,i_t,j_t)$. Then $(i_t,j_t)$ is publicly observed and the play goes to stage $t+1$. Notice that the payoff $g(k_t,l_t,i_t,j_t)$ is not directly observed and may not be deduced by the players at the end of stage $t$.

We denote by $\Sigma$ and ${\cal T}$ the sets of behavior strategies of the players. Formally, for $t \geq 1$, let  $H_t=(I\times J)^t$  denote the set of histories of actions of length $t$, with the convention $H_0=\{\emptyset\}$. A strategy $\sigma \in \Sigma$ is a sequence $(\sigma_t)_{t \geq 1}$, where $\sigma_t$ is a map from $H_{t-1}\times K^t$ to $\Delta(I)$. For $T \in \N^*$, the set $\Sigma_T$ denotes the set of $T$-stage strategies for Player 1, that is of finite sequences $(\sigma_1,..,\sigma_T)$ induced by elements of $\Sigma$. $\mathcal{T}$ and $\mathcal{T}_T$ are defined similarly.

The value of the $T$-stage game is denoted $v_T(p,q)$, i.e. 
\[v_T(p,q)= \max_{\sigma \in  \Sigma}\min_{\tau \in {\cal T}} \gamma_T^{p,q}(\sigma, \tau)=\min_{\tau \in {\cal T}} \max_{\sigma \in  \Sigma}\gamma_T^{p,q}(\sigma, \tau),\]
where $\gamma_T^{p,q}(\sigma, \tau)=\E_{p,q, \sigma, \tau} [\frac{1}{ T} \sum_{t=1}^T g(k_t,l_t,i_t,j_t)]$, the expectation being taken with respect to the probability $\P_{p,q,\sigma,\tau}$ induced by $(p,q,\sigma,\tau)$ on the set of plays $\Omega=(K\times L \times I\times J)^\infty$.  Elements of $\Delta(K)$  are seen as row vectors, so that if $p'$ in $\Delta(K)$  is the law of the state $k_T$ of stage $T$, then $p'M$ is the law of the  following state $k_{T+1}$ (similarly elements of $\Delta(L)$ are seen as row vectors). We have the standard recursive formula, for all $T\geq 0$: 
\begin{align*}
(T+1)v_{T+1}(p,q) &=\max_{x \in \Delta(I)^K} \min_{y \in \Delta(J)^L} \left(G(p,q,x,y) + T \sum_{i\in I, j \in J} x(p)(i)y(q)(j) v_{T}( \vec{p}(x,i)M, \vec{q}(y,j)N)\right)\\
&=\min_{y \in \Delta(J)^L} \max_{x \in \Delta(I)^K}  \left(G(p,q,x,y) + T \sum_{i\in I, j \in J} x(p)(i)y(q)(j) v_{T}( \vec{p}(x,i)M, \vec{q}(y,j)N)\right).
\end{align*}
where $G(p,q,x,y)=\sum_{k,l,i,j} p^k q^l x^k(i)y^l(j) g(k,l,i,j)$,  
\begin{align*} 
\forall i \in I, \;x(p)(i)=\sum_k p^k x^k(i), \; \text{and} \;  \vec{p}(x,i)= \left(\frac{p^kx^k(i) }{x(p)(i)} \right)_{k\in K} \in \Delta(K) \;\text{if} \; x(p)(i)>0, \\
\forall j \in J, \; y(q)(j)=\sum_l q^l y^l(j), \;\text{and} \; \vec{q}(y,j)=\left(\frac{q^ly^l(j) }{y(q)(j)} \right)_{l\in L} \in \Delta(L) \;\text{if} \; y(q)(j)>0,
\end{align*}
and $\vec{p}(x,i)$ (resp. $\vec{q}(y,j)$)  is defined arbitrarily if $x(p)(i)=0$ (resp. $y(q)(j)=0$).    If at some stage $t$  the belief of Player 2 on the current state $k_t$ in $K$ is represented by $p$, and if Player 1 plays    the mixed action $x^k$ if his  current state is $k$, then $x(p)(i)$ will be the probability that a given  action $i$ in $I$ is played at this stage. After observing $i$, $\vec{p}(x,i)$   represents the new belief of Player 2 on $k_t$, and  $\vec{p}(x,i)M$ is the belief  of Player 2 on $k_{t+1}$.\\

We will extensively use the following definitions. If $f$ is a real function on   $\Delta(K)\times \Delta(L)$, we say that $f$  is  $I$-concave, resp. $\two$-convex,  if it is concave in the first variable, resp. convex in the second variable. $\cav_I f$ and $\vex_{\two}f$  denote respectively the smallest $I$-concave function above $f$ and the largest $\two$-convex function below $f$. 
\begin{defi}   ${\cal C}$  denotes the set of continuous functions from $\Delta(K)\times \Delta(L)$ to $[-1,1]$. For  each $f$ in ${\cal C}$, we introduce: 

\noindent $\mathcal{C}^+(f) = \left\{   w \in {\cal C} , w\; \text{is $I$-concave \; and\; } 
   \forall (p,q)\in \Delta(K)\times \Delta(L),\; w(p,q) \geq \vex_{\two} \Max(w,f)(p,q)  \right\}, $
\noindent    $\mathcal{C}^-(f) = \left\{   w \in {\cal C} , w\; \text{is $\two$-convex \; and\; } 
   \forall (p,q)\in \Delta(K)\times \Delta(L),\; w(p,q) \leq \cav_I \Min(w,f)(p,q)  \right\}. $\end{defi}
 
  \vspace{0,3cm}
  
The main result of the paper is the following. \\
  \begin{thm} \label{thm1} $v(p,q)=\lim_{T \rightarrow \infty} v_T(p,q)$ exists. \end{thm}
 
 \vspace{0,3cm}
 
Moreover,  if   the Markov chains on $K$ and $L$ are recurrent and aperiodic, we will define   a ``nonrevealing limit value" function $\hat{v}$ (see Proposition \ref{prohatv}), and will prove that $v$ is  the unique solution of the Mertens-Zamir system associated to $\hat{v}$, that is:
 \begin{equation}\label{systemeMZ} 
\forall (p,q) \in \Delta(K)\times \Delta(L),   \left\{ \begin{matrix} w(p,q) &=& \vex_{\two} \Max(w,\hat{v})(p,q) \\
w(p,q) &=& \cav_{I} \Min(w,\hat{v})(p,q) \end{matrix} \right. ,
\end{equation}
\noindent or, equivalently, that  $v   = \inf\{  w \in {\cal C}^+(\hat{v})\}    = \sup   \{  w \in  {\cal C}^-(\hat{v})\}.$
\p
Note that this characterization is helpful only when one of  the chains has several recurrence classes, i.e. when there is incomplete information on the asymptotic behavior of the chains. Indeed, in the case of irreducible chains, our definition of $\hat{v}$ implies that $v(p,q)=\hat{v}(p,q)$ and it is easy to show that these functions do not depend on $(p,q)$ (see example C in the next section), so that the above system do not bring any information.

\section{Preliminaries}    \label{sec3}
 In all the paper the sets $\Delta(K)$, $\Delta(L)$ and $\Delta(K)\times \Delta(L)$ are endowed with the $L^1$-norm.
 The following properties of the value functions are standard. For each $T\geq 1$, $v_T$ is $I$-concave, $\two$-convex and 1-Lipschitz. And for all   $(p,q) \in \Delta(K)\times \Delta(L)$, we have:
$|v_{T+1}(p,q)-v_{T}(p,q)| \leq \frac{2}{T}, \;\; |v_{T+1}(p,q)-v_{T}(pM,q N)| \leq \frac{2}{T}, $ 
and  consequently $|v_{T}(p,q)-v_{T}(pM,qN)| \leq \frac{4}{T}$. \\

If some state $k$ is transient, then for each $\varepsilon$ there exists $T_0$ such that $\P ( \forall T \geq T_0,\; k_T\neq k)\geq 1-\varepsilon$. Both  players can wait until stage $T_0$ and enter a game where state $k$ has disappeared  with high probability, so to prove theorem \ref{thm1} we can assume w.l.o.g. that there is no transient state. In the sequel, we assume that all states in $K$ and $L$ are recurrent.  \\
 
We now take care of periodicity.  Some states may {\it a priori} be periodic, let $T_0$ be a common multiple of all periods in the chains induced by $M$ and $N$.  From the above properties, we   have $\limsup_{T}v_{TT_0}(p,q)=\limsup_{T}v_{T}(p,q)$, and $\liminf_{T}v_{TT_0}(p,q)=\liminf_{T}v_{T}(p,q)$, so to study $\lim_T v_T$ we can consider plays by blocks of $T_0$ stages (so that one stage now corresponds to $T_0$ original stages). The Markov chains with transitions matrices $M^{T_0}$ and $N^{T_0}$ are   aperiodic, and this consideration is without loss of generality.\\

Summing up, from now on\footnote{except in the example  E   of the next section  and in Remark \ref {cavtransient}}  we assume w.l.o.g. that both Markov chains are recurrent and aperiodic. We still use the letters $M$ and $N$ to denote the (now recurrent aperiodic) transitions matrices. We have the convergence of $M^t$ to a stochastic matrix $B$ such that $B^2=B=BM=MB$. The state space $K$ can be partitioned into recurrence classes $K(1)$, ..., $K(r_M)$, where $r_M$ is the number of recurrence classes associated to $M$. Each class  $K(r)$ has a unique invariant measure $p^*(r)$, and the $k$-th row of $B$ corresponds to $p^*(r)$ where $K(r)$ is the  class containing $k$. Invariant measures for $B$ and $M$ coincide, the set of those measures is   the convex hull of $\{p^*(1),...,p^*(r_M)\}$ denoted by $P^*$. Although $B$ is not invertible if there are less recurrence classes than states, for $p^*$ in $P^*$  we write $B^{-1}(p^*)=\{p\in \Delta(K), pB=p^*\}.$ 

Regarding the Markov chain associated to the states in $L$, we proceed similarly and use the following notations. There is convergence of $N^t$ to a stochastic matrix $C$. The state space $L$ is partitioned into recurrence classes $L(1)$,  ..., $L(r_L)$. Each class  $L(r)$ has a unique invariant measure $q^*(r)$, and the $l$-th row of $C$ corresponds to $q^*(r)$ where  $L(r)$ is the  class containing $l$.  Invariant measures for $C$ and $N$ coincide, the set of those measures is   the convex hull of $\{q^*(1),...,q^*(r_L)\}$ denoted by $Q^*$. For $q^*$ in $Q^*$, we write $C^{-1}(q^*)=\{q\in \Delta(K), qC=q^*\}.$ \\

We now provide an example to  illustrate    some of our notations.

\noindent {\bf Example A:}   $K=L=\{a,b,c\}$, and 
$ M=N=\begin{pmatrix}  \frac{2}{3} & \frac{1}{3} & 0 \\ \frac{1}{3} & \frac{2}{3} & 0 \\ 0
&0 & 1 \end{pmatrix} .$ There are two recurrence classes $K(1)=\{a,b\}$ and $K(2)=\{c\}$. The associated invariant measures are $p^*(1)=(
\frac{1}{2},\frac{1}{2},0)$ and $p^*(2)=(0,0,1)$, and the limit projection matrix is   
$B= \begin{pmatrix} \frac{1}{2} &\frac{1}{2}& 0 \\ \frac{1}{2} & \frac{1}{2} &0 \\ 0 &0 &1 \end{pmatrix}. $
 
\begin{center}
\includegraphics[scale=0.5]{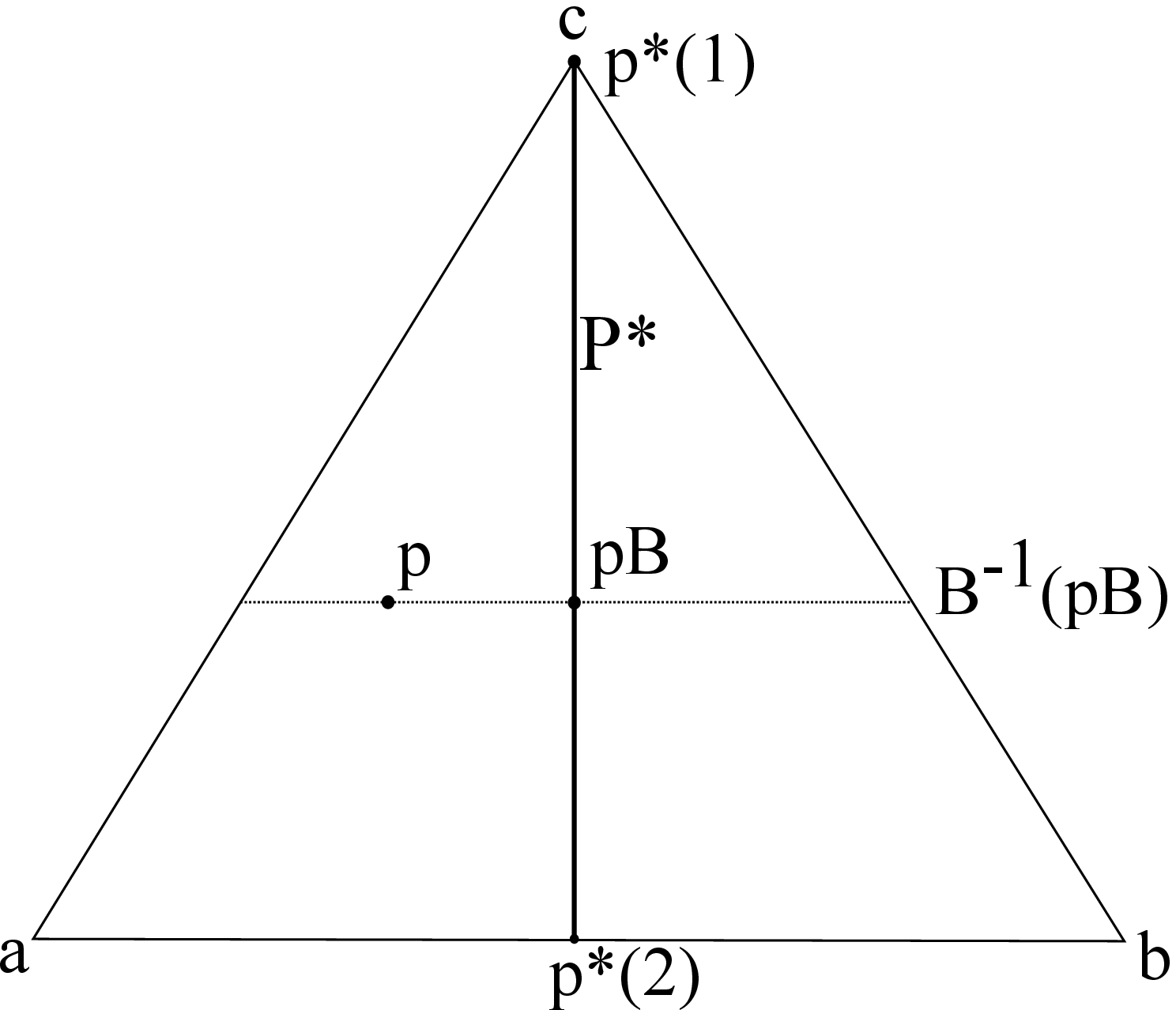}
\end{center}

\vspace{0,3cm}

We will use the following notion in the proofs.  \begin{defi}  $f$ is 
 {balanced}  if for all $(p,q) \in  \Delta(K)\times \Delta(L)$,  $f(p,q)=f(pM,qN).$ \end{defi}
 Since the Markov chains are recurrent aperiodic, $f$ is balanced if and only if $f$ is constant on the recurrence classes of the chains.

\section{Examples  and overview of the proof} \label{sec4}

Let us at first emphasize that, in contrast with repeated games with incomplete information, if Player 1 does not use his information on the states $(k_t)_{t \geq 1}$, beliefs of Player 2 over states evolve across stages according to the transition matrix $M$. 
Precisely, under this assumption, the belief of Player 2 on the current state $k_{t+1}$ at stage $t + 1$ (i.e., the conditional law of $k_{t+1}$ given the actions played up to stage $t$) is simply the law of  $k_{t+1}$, which is $pM^{t}$. Symmetrically, if Player 2 does not use his information on the states $(l_t)_{t \geq 1}$, then the belief of Player 1 on $l_{t+1}$ at stage $t+1$ is $qN^{t}$.\\

\noindent {\bf Example B: A Markov chain without memory.} If both Markov chains are i.i.d. sequences, then the problem reduces to a classical repeated game, by replacing the payoff function in pure strategies by its expectation with respect to $(k_t,l_t)$. The value of the infinitely repeated game exists and equals the value of this auxiliary one-stage game.\\

\noindent {\bf Example C:  Irreducible aperiodic Markov chains.} 
Assume that both chains are irreducible and aperiodic, say 
$ K=L=\{a,b\}$, $M = \begin{pmatrix} \frac{2}{3} &\frac{1}{3} \\ \frac{1}{3} & \frac{2}{3} \end{pmatrix}$ and $  N = \begin{pmatrix} \frac{3}{4} &\frac{1}{4} \\ 1 & 0 \end{pmatrix}. $
This Markov chain has a unique recurrence class, so $r_M=r_N=1$ and $K=K(1)=L=L(1)$. Unique invariant probabilities  are  $p^*=( \frac{1}{2},\frac{1}{2})$  and $q^*=( \frac{4}{5},\frac{1}{5})$.   $M^t$ and $N^t$  respectively converge to:
$B= \begin{pmatrix} \frac{1}{2} & \frac{1}{2} \\ \frac{1}{2} & \frac{1}{2} \end{pmatrix}$ and $ C= \begin{pmatrix} \frac{4}{5} & \frac{1}{5} \\ \frac{4}{5} & \frac{1}{5} \end{pmatrix}. $
In this case, the limit value exists, and we even have the stronger result that the (uniform) value of the infinitely repeated game exists. At first, if Player 1 plays independently of the states, the beliefs of Player 2 converge to $p^*$. 
On the other hand, after any sufficiently long number of stages, he can always forget the past actions of Player 2 and act as if his beliefs were very close to $q^*$. Choose now some integer $T_0$ such that $M^{T_0}$ and $N^{T_0}$ are respectively close to $B$ and $C$, and then $T_1$ much larger than $T_0$. Let us construct a strategy of Player 1 as follows: (1) play any strategy independent of the states during  $T_0$ stages, then (2) play an optimal strategy in the game $\Gamma_{T_1}(p^*,q^*)$, and (3) come back to (1). Such a strategy guarantees to Player 1 $v_{T_1}(p^*,q^*)$ up to some arbitrarily small error, hence by choosing $T_1$ large enough Player 1 guarantees $\lims_T v_T(p^*,q^*)$ in the infinitely repeated game. Inverting the role of the players implies that the value of the infinitely repeated game  exists and equals  $\lim_T v_T(p^*,q^*)$. In particular  it  does not depend on the initial probabilities $(p,q)$.  Notice that this proof can be generalized to the case of any irreducible aperiodic Markov chain. \\

 \noindent {\bf Example D:  $M$ and $L$  are  the identity matrix.} In this case, the problem is a repeated game with incomplete information on both sides as studied by Mertens and Zamir \cite{mertenszamir} since the initial state is selected at the beginning of the game and remains constant. In terms of Markov chains, the states\footnote{also called {\it types} in  games with incomplete information}  can be here identified with the recurrence classes, and therefore any information on the states is an information on the asymptotic behavior of the corresponding chain.

Let us at first define the non-revealing game:  it denotes the game where none of the players is allowed to use his information on the state variable (i.e. players are restricted to use strategies which do not depend on the selected states).  For each pair $(p,q)$ of initial probabilities, this game can be analyzed as the repetition of the average matrix game $\sum_{k,l} p^k q^l g(k,l,.,.)$. 
Consequently, its value is just the value $u(p,q)= \max_{x \in \Delta(I)} \min_{y \in \Delta(J)}  \sum_{k,l} p^k q^l g(k,l,x(i),y(j))$ of the one-stage average game.

\begin{thm}[\cite{mertenszamir}]\label{thmmz} 
For each $f$ in ${\cal C}$, there exists a unique solution in ${\cal C}$   to the system: 
 \begin{equation}\label{systemeMZ2} 
\forall (p,q) \in \Delta(K)\times \Delta(L),   \left\{ \begin{matrix} w  &=& \vex_{\two} \Max(w,f)  \\
w  &=& \cav_{I} \Min(w,f) \end{matrix} \right. ,
\end{equation}

  This solution is denoted $MZ(f)$ and $MZ(f)= \sup \{ w | w \in \mathcal{C}^-(f) \} = \inf \{ w | w \in \mathcal{C}^+(f)\}$.
\end{thm}

\vspace{0,3cm}

Back to games, Mertens and Zamir proved that in the case of example $D$, the limit value $v$ exists and is  $MZ(u)$.      Assuming the equality  $\sup \{ w | w \in \mathcal{C}^-(u) \} = \inf \{ w | w \in \mathcal{C}^+(u)\}$, their proof of convergence  of $(v_T)$  can be sketched as follows. The main point is to show that for any function $w$ in $\mathcal{C}^-(u)$, Player 1 can defend   in the game of length $T$ with initial probabilities  $p$ and $q$, the quantity $w(p,q)$  up to an error going to zero with $T$ uniformly with respect to $(p,q)$ and the chosen strategy of Player 2. This can be done as follows. The strategy of Player 2 being known, Player 1 can compute after each stage $t$  the current pair of beliefs $(p_t,q_t)$ of both players: $p_t$ denotes the conditional law of $k_{t+1}$ given past actions and $q_t$ denotes the conditional law of $l_{t+1}$ given past actions. Then, either $u(p_t,q_t) \geq w(p_t,q_t)$ and he plays an optimal strategy in the non-revealing game at $(p_t,q_t)$, or $u(p_t,q_t) < w(p_t,q_t)$ and he can split his information  (think of Player 1 sending a, possibly random, message to Player 2)  in order to drive Player 2's  beliefs at some points $\tilde{p}_t$ such that $u(\tilde{p}_t,q_t) \geq w(\tilde{p}_t,q_t)$. The existence of such a splitting is ensured by: 1)  the property $w\in \mathcal{C}^-(u)$  and   2)   the splitting Lemma (see e.g. the Lemmas  \ref{splitting} and \ref{SLemma2} in the present work). This strategy would defend $w(p_t,q_t)$ at stage $t+1$ if Player 2 was not using his information on the state at stage $t+1$. This is almost the case since the error can be bounded by the expected $L_1$-variation of the beliefs $\|q_{t+1}-q_t\|$. Recall that the sequence of beliefs $(q_t)$ is a martingale and a classical bound on its $L_1$ variation allows to bound uniformly the average error by a quantity vanishing with $T$. Moreover  using 
the properties of $w$ and the  construction of the strategy of Player 1,    the expected value of $w(p_t,q_t)$ is always greater or equal than $w(p,q)$. All this implies that  $\limi_n v_n(p,q) \geq w(p,q)$, so we obtain    $\limi_n v_n \geq \sup\{ w | w \in \mathcal{C}^-(u)\}$. 
Inverting the role of the players, we also have $\lims_n v_n \leq \inf\{ w | w \in \mathcal{C}^+(u)\}$ which concludes the proof.  \\

\noindent {\bf Example E: A periodic chain}  Let $K$ be $\{a,b\}$, and $M$ be  $\left(\begin{array}{cc}0&1\\1&0\end{array}\right)$.

This Markov chain has a unique recurrence class, which is periodic with period $2$. 
The relevant information for Player 2 is not the recurrence class, but the initial state which determines whether the sequence of states $(k_t)_{t \geq 1}$ will be $\{a,b,a,b,...\}$ or $\{b,a,b,a,...\}$. By considering the auxiliary game in which stages are blocks of length $2$ of the initial game, then the new transition matrix is $M^2$, hence the identity, and each initial state becomes now a recurrence class. The problem can be reduced   to that of Example D.  As explained in  section \ref{sec3}, this idea of playing stages by blocks of fixed length can be used in the general case to assume without
loss of generality that the chain is aperiodic. \\

\noindent Back to {\bf  Example A}, the existence of the limit value in this case    appears to be more difficult, and we do not know how to proceed
without following the general proof presented in this paper. Following the general scheme of Renault \cite{R2006}, we will consider two kinds of information: long-term information (corresponding to the recurrence classes of the two chains, as in the Mertens-Zamir case),
and short-term information (corresponding to the state variable itself within some recurrence class). \\

We now come back to the general case, assuming the Markov chains are recurrent and aperiodic,  and describe the different steps of our proof.

\subsection*{Overview of the Proof}

By analogy with the  proof in Renault \cite{R2006}, we define a non-revealing game, where Players are restricted to strategies that do not reveal information on the recurrence classes of the chains. Precisely, the strategy of Player 1 has to be such that, after each history of actions of length $t$, the conditional probability that $k_{t+1}$ belongs to any recurrence class $K(r)$ remains unchanged (see Definition \ref{stratNR}). 

This conceptual definition has to be compared with the non-revealing strategies introduced by Aumann and Maschler in case each player observes after each stage a signal which is a function of the actions played, and which requires that the beliefs of Player 2 on the state $k$ remains unchanged. 

The main difference here is that the study of the non-revealing game can not be reduced to that of a finite matrix game. 
Following the scheme of the proof of example C, with appropriate modifications, we show that the limit value $\hat{v}$ of this non-revealing game exists and that it depends only on the limit distributions of the two Markov chains $(pB,qC)$.
This proof requires however a precise analysis of the $T$-stage  non-revealing games, for which we prove the existence of the value $\hat{v}_T$ and an appropriate  recursive formula. The crucial point of the analysis of the non revealing values is to establish that the family  $(\hat{v}_T)_{T \geq 1}$  is uniformly equi-continuous (this is the aim of the whole section \ref{proofvhat}), which implies uniform convergence of these values to their limit. 

The second main technical  difficulty is to link the values of the non revealing games with those of the original game, i.e. to adapt the    proof of Mertens and Zamir illustrated  in example D.  We use  the uniform convergence  mentioned above, in particular   for sufficiently large $T_0$, $\hat{v}_{T_0}$ is $\varepsilon$-close to $\hat{v}$. We prove that for any $n \in \N^*$, and for any {\it balanced} function $w \in \mathcal{C}^-(\hat{v})$, Player 1 can defend $w(p,q)$ in $\Gamma_{nT_0}(p,q)$  up to some error going to zero with $n$ uniformly with respect to $(p,q)$ and the chosen strategy of Player 2. The construction of the strategy of Player 1 is  the same as explained in example D, except that stages are replaced by blocks of $T_0$ stages: Player 1 can compute after each block the current pair of beliefs $(p_{nT_0},q_{nT_0})$ of both players and either (1) $\hat{v}(p_{nT_0},q_{nT_0}) \geq w(p_{nT_0},q_{nT_0})$ and he plays an optimal strategy in the non-revealing game of length $T_0$ at $(p_{nT_0},q_{nT_0})$, or (2) $\hat{v}(p_{nT_0},q_{nT_0}) < w(p_{nT_0},q_{nT_0})$, and  he can send a random  message to Player 2 in order to drive his beliefs at some points $\tilde{p}_{nT_0}$ such that $\hat{v}(\tilde{p}_{nT_0},q_{nT_0}) \geq w(\tilde{p}_{nT_0},q_{nT_0})$ and there  he plays as in case (1). This strategy would defend $w(p_{nT_0},q_{nT_0})$ up to an error of $\varepsilon$ on the $n$-th block of size $T_0$ if Player 2 was not using his information on the recurrence classes of the chain $(l_t)_t$ on this block. In order to bound the error due to this approximation, we have to replace the error term based on the $L_1$-variation of the process $(q_t)_{t \geq 1}$ (which is no more a martingale here)  by the $L_1$-variation of an auxiliary martingale representing the beliefs of Player 1 over the recurrence classes of the Markov chain $(l_t)_{t\geq 1}$ (note however that we obtain a less precise bound since the method developed for games with incomplete information does not apply here). Then, we prove that  the expected value of $w(p_{nT_0},q_{nT_0})$ is always greater or equal to $w(p,q)$ using similar tools as in \cite{mertenszamir} and also that  $w$ is balanced. We obtain $\limi_n v_n \geq \sup\{ w | w \in \mathcal{C}^-(\hat{v}), {\rm balanced}\}$, and by symmetry   $\lims_n v_n \leq \inf\{ w | w \in \mathcal{C}^+(\hat{v}), {\rm balanced}\}$. 

Then it remains to show  the equality 
\[   \sup\{ w | w \in \mathcal{C}^-(\hat{v}), {\rm balanced} \}= \inf \{w | w \in \mathcal{C}^+(\hat{v}), {\rm balanced} \}, \] 
by reducing the problem to the equality showed by Mertens and Zamir. The convergence of $(v_n)_n$  follows, and it is not difficult to see that the limit is $MZ(\hat{v})$.

\section{Non revealing games}\label{nonrevealing}

Let $\sigma$ in $\Sigma$ be a strategy of Player 1, i.e. $\sigma=(\sigma_T)_{T\geq 1}$, with $\sigma_T: (K \times I \times J)^{T-1}\times K\longrightarrow \Delta(I)$ for all $T$. For each finite history $h_T$ in $(I \times J)^{T}$, we say that  $h_T$ is compatible with $(p, \sigma)$ if for some $\tau$ and $q$ we have $\P_{p,q,\sigma, \tau}(h_T)>0$. 
In this case we define $p_{T}(p, \sigma)(h_T)$ in $\Delta(K)$ as the law of the state $k_{T+1}$ knowing that: the initial state $k_1$ is selected according to $p$, Player 1 uses $\sigma$ and $h_T$ has been played at the first $T$ stages.  It does not depend  on the last move of Player 2 in $h_T$, and for all $k$ in $K$, $q$ in $\Delta(L)$ and strategy  $\tau$ of Player 2 such that $\P_{p,q,\sigma, \tau}(h_T)>0$, we have:
\[ p^k_{T}(p, \sigma)(h_T)=\P_{p,q,\sigma, \tau}(k_{T+1}=k|h_T). \]
If $h_T$ is not compatible with $(p,\sigma)$, we define $p_T(p,\sigma)(h_T)$   arbitrarily in $\Delta(K)$.   For $T=0$,  $p_0(p, \sigma)=p$ is the law of the first state $k_1$. 

\vspace{0,3cm}

\begin{defi} We define $\hat{p}_{T}(p, \sigma)(h_T)=p_{T}(p, \sigma)(h_T)B$, and we call $\hat{p}_{T}(p, \sigma)(h_T)$ the relevant information of Player 2 after $h_T$ has been played. \end{defi}

\vspace{0,3cm}
 $\hat{p}_{T}(p, \sigma)(h_T)$ indicates  the belief on the asymptotic distribution of the Markov chain $(k_t)$  after $h_T$ has been played, if the initial probability is $p$ and player 1 uses $\sigma$.  
 
For each $T$, we denote by  ${\cal H}_T$   the $\sigma$-algebra on the set of plays $\Omega$ generated by the projection on $(I\times J)^T$ giving the first $T$ moves.
  
\begin{lem}
For any strategy pair $(\sigma, \tau)$ in $\Sigma \times  {\cal T}$, the process $(\hat{p}_T(p, \sigma))_{T\geq 0}$ is a $({\cal H}_T)_{T\geq 0}$-martingale with respect to $\P_{p,q,\sigma, \tau}.$ 
\end{lem}

\vspace{0,3cm}

\noindent{Proof:} 
Assume $(\sigma,\tau)$ is played, the    initial  probabilities  on $\Delta(K)$ and $\Delta(L)$ being $p$ and $q$. Fix $T\geq 0$  and $h_T=(i_1,j_1,...,i_T,j_T)\in (I\times J)^T$. Given  $(i_{T+1},j_{T+1})$  in $I\times J$, we  write $h_{T+1}=(i_1,...,j_T,i_{T+1},j_{T+1})$ and denote by  $r_T(p,\sigma)(h_{T+1})$    the conditional probability on the state $k_{T+1}$ given that   $h_{T+1}$ has been  played at the first $T+1$ stages. We have $r_T(p,\sigma)(h_{T+1}) M= p_{T+1}(p,\sigma)(h_{T+1})$, and 
$\E_{p,q,\sigma,\tau}(r_T(p,\sigma)|h_T)=p_T(p,\sigma)(h_T)$. Now,
\begin{align*}
\E_{p,q,\sigma,\tau}(\hat{p}_{T+1}(p,\sigma)|h_T) & = \E_{p,q,\sigma,\tau}({p}_{T+1}(p,\sigma)B|h_T) \\
 & = \E_{p,q,\sigma,\tau}({r}_{T}(p,\sigma)MB|h_T) \\
  & = \E_{p,q,\sigma,\tau}({r}_{T}(p,\sigma)B|h_T) \\
  & = p_{T}(p,\sigma)(h_T) B\\
   & = \hat{p}_{T}(p,\sigma)(h_T). 
\end{align*}

\begin{defi}\label{stratNR}
 A strategy $\sigma$ in $\Sigma$  is called non revealing at $p$ if for all $\tau$ and $q$, the martingale $(\hat{p}_T(p, \sigma))_{T\geq 0}$ is $\P_{p, q, \sigma, \tau}$-almost surely constant. The set of NR strategies of Player 1 at $p$ is denoted $\hat{\Sigma}(p)$.
\end{defi}

We will  give an alternative definition of non revealing strategies. 

\begin{defi}
\begin{align*}
NR(p) & = \{x \in \Delta(I)^K, \forall i \in I \mbox{ s.t. } x(p)(i)>0, \;\; \vec{p}(x,i)B =pB \} \\
& = \{x \in \Delta(I)^K, \forall i \in I, (p^kx^k(i))_{k\in K} B =x(p)(i)\, pB \}.
\end{align*}\end{defi}
Given $\sigma$ in $\Sigma$ and  an history $h_T\in (I\times J)^T$ which is  compatible with $(p, \sigma)$, for each  $i$ in $I$ and $k$ in $K$, we denote by  $x^k(p,\sigma)(h_T)(i)$ the probability $\P_{p,q,\sigma,\tau} (i_{T+1}=i | k_{T+1}=k,h_T)$ that Player 1 plays $i$ after  $h_T$ if the next state is $k$.  This probability   does not depend on  $(q,\tau)$. The vector $x (p,\sigma)(h_T)= (x^k(p,\sigma)(h_T)(i))_{k,i}$ is viewed as en element of $\Delta(I)^K$. The  proof of the next lemma is similar to the proof of Proposition 6.5 in \cite{R2006} and is omitted here. \\

\begin{lem} \label{NRSlem}
A strategy $\sigma$ in $\Sigma$ is non revealing at $p$ if and only if for all $T\geq 0$ and history $h_T\in (I\times J)^T$ compatible with $(p, \sigma)$, the vector    $x (p,\sigma)(h_T)$ belongs to $NR(p_T(p,\sigma)(h_T))$.
\end{lem}

  We  also have a splitting lemma with non revealing strategies.   
\begin{lem} \label{SLemma}
Consider  a convex combination $p=\sum_{s \in S}\alpha_s p_s$ in $\Delta(K)$ having the property that $p_sB=pB$ for all $s$.
Consider for each $s$ a non revealing strategy $\sigma_s$ in $\hat{\Sigma}(p_s)$. There exists $\sigma$ in $\hat{\Sigma}(p)$ such that: 
\[\forall q \in \Delta(L), \forall \tau \in {\cal T},\;\; \P_{p,q,\sigma,\tau}=\sum_{s \in S}  \alpha_s \P_{p_s,q,\sigma_s,\tau}.\]
\end{lem}

 The proof is similar to the proof of Lemma 6.6 in \cite{R2006}, with $\sigma$ being defined via the splitting procedure: observe the first  state $k_1$ in $K$, then choose $s$ according  to the probability $\alpha_s p_s^{k_1}/p^{k_1}$, and play according  to $\sigma_s$. 
 \\
Non revealing strategies of Player 2 at $q$  are defined similarly, the set of such strategies is denoted  $\hat{\cal T}(q)$.
\begin{defi}
The $T$-stage non revealing game at $(p,q)$ is the $T$-stage game where the strategy sets are restricted to $\hat{\Sigma}(p)$ and $\hat{\cal T}(q)$. It is denoted by $\hat{\Gamma}_T(p,q)$. 
\end{defi}
 Define $\hat{\Sigma}_T(p)$ as the set of $T$-stage non revealing strategies of Player 1 at $p$. Formally, $\hat{\Sigma}_T(p)$ is the projection of $\hat{\Sigma}(p)$ over the (compact) set of $T$-stage strategies of Player 1. Similarly, $\hat{\cal T}_T(q)$ denotes the set of $T$-stage non revealing strategies of Player 2 at $q$. The game  $\hat{\Gamma}_T(p,q)$ can equivalently be seen as the zero-sum game  $(\hat{\Sigma}_T(p), \hat{\cal T}_T(q), \gamma_T^{p,q})$. 
The proof of the next proposition is standard and similar to the proofs of proposition 7.3 and 7.4 in  \cite{R2006}.
\begin{pro} \label{pro1}For all $T\geq 1$, the $T$-stage non revealing game at $(p,q)$ has a value denoted by $\hat{v}_T(p,q)$. 
 For each $q$ in $\Delta(L)$ and $p^*$ in $P^*$, $(p\mapsto \hat{v}_T(p,q))$ is concave on $B^{-1}(p^*)$, and similarly for each $p$ in $\Delta(K)$ and $q^*$ in $Q^*$, $(q\mapsto \hat{v}_T(p,q))$ is convex on $C^{-1}(q^*). $ Moreover,  for all $T\geq 0$ we have the recursive formula: 
\begin{align*}
(T+1)\hat{v}_{T+1}(p,q) &=\max_{x \in NR(p)} \min_{y \in NR(q)} \left(G(p,q,x,y) + T \sum_{i\in I, j \in J} x(p)(i)y(q)(j) \hat{v}_{T}( \vec{p}(x,i)M, \vec{q}(y,j)N)\right)\\
&=\min_{y \in NR(q)} \max_{x \in NR(p)}  \left(G(p,q,x,y) + T \sum_{i\in I, j \in J} x(p)(i)y(q)(j) \hat{v}_{T}( \vec{p}(x,i)M, \vec{q}(y,j)N)\right).
\end{align*} 
\end{pro}

We  now show that   to play optimally,  Player 2 does not need to remember the whole sequence of states $(l_1,...,l_t,...)$. This will allow us   to consider a  smaller    set  of strategies. \\
\begin{defi} 
Let $\overline{\mathcal{T}}$ denote the subset of strategies of Player $2$ which for all $T\geq 0$, depend at stage $T+1$ only on the past history $h_T \in (I \times J)^T$ and on the current state $l_{T+1}$. Define $\overline{\Sigma}$ similarly.
\end{defi}
\begin{pro} \label{pro27}
For all $(p,q) \in \Delta(K) \times \Delta(L)$, for all $\tau \in \mathcal{T}$, there exists a strategy $\overline{\tau} \in  \overline{\mathcal{T}}$
such that for all $T \geq 1$ and $\sigma \in \Sigma$, $\gamma_T^{p,q}(\sigma,\tau)=\gamma_T^{p,q}(\sigma,\overline{\tau})$.
Moreover, if $\tau$ is non-revealing at $q$, then $\overline{\tau}$ is also non-revealing at $q$. As a corollary, in $\Gamma_T(p,q)$ and in the non revealing game $\hat{\Gamma}_T(p,q)$  the players have optimal strategies in $\overline{\Sigma}$ and $\overline{\mathcal{T}}$.
\end{pro}

\vspace{0,3cm}

\noindent {Proof:}
The strategy $\overline{\tau}$ proceeds as follows:  for $t\geq0$, at stage $t+1$, Player 2 does not remind the past states $(l_1,..,l_{t})$, but using the history $h_t$, he can generate a virtual sequence $(\overline{l}_1,..,\overline{l}_t)$ using the conditional law of $(l_1,..,l_t)$ given $(h_t,l_{t+1})$ under $\P_{p,q,\sigma,\tau}$ (which does not depend on $\sigma$). He selects then at stage $t+1$ an action $j_{t+1}$ with the probability $\tau(h_t,\overline{l}_1,..,\overline{l}_t,l_{t+1})$. 
Formally,  we have $\overline{\tau}(h_{t})(l_{t+1}) = \E_{p,q,\sigma,\tau}[ \tau(h_t,l_1,...,l_t,l_{t+1}) | h_{t},l_{t+1}].$ 

We will prove that for all $t\geq 1$, $(h_t,l_t,k_t)$ has the same distribution under $\P_{p,q,\sigma,\tau}$ and under $\P_{p,q,\sigma , \overline{\tau}}$, i.e. for all $h \in H_{t}$, for all $(k,l) \in K\times L$, 
\begin{equation}\label{induction}
\P_{p,q,\sigma,\tau}(h_{t}=h,l_{t}=l,k_{t}=k)=\P_{p,q,\sigma,\overline{\tau}}(h_{t}=h,l_{t}=l,k_{t}=k). 
\end{equation}
Let us proceed by induction on $t$. The property is obvious for $t=1$, 
assume it is true  for some $t \geq 1$.  
At first, the conditional distribution of $l_{t+1}$ given $(h_{t},k_{t},l_{t})$ is $M(l_{t},.)$ under $\P_{p,q,\sigma,\tau}$ and under $\P_{p,q,\sigma , \overline{\tau}}$. We deduce that the law of $(h_t,l_{t+1})$ is the same under both probabilities. By construction the conditional law of $j_{t+1}$ given $(h_t,l_{t+1})$ is the same under $\P_{p,q,\sigma,\tau}$ and under $\P_{p,q,\sigma , \overline{\tau}}$, which implies that the conditional law of $(l_{t+1},j_{t+1})$ given $h_t$ is also the same. Using that $(k_{t+1},i_{t+1})$ and $(l_{t+1},j_{t+1})$ are  conditionally independent given $h_t$ (under both distributions) and that the conditional law of  $(k_{t+1},i_{t+1})$ given $h_t$ does not depend on $\tau$, the conditional distribution of $(k_{t+1},l_{t+1,},i_{t+1},j_{t+1})$ given $h_t$ is the same under $\P_{p,q,\sigma,\tau}$ and under $\P_{p,q,\sigma , \overline{\tau}}$. We deduce that the law of $(h_{t+1},k_{t+1},l_{t+1})$ is the same under both probability distributions which concludes the proof of \eqref{induction}.
Going back to the main proof, we obtain
\begin{align*} 
\gamma_T^{p,q}(\sigma,\tau) &= \sum_{t=0}^{T-1} \E_{p,q,\sigma,\tau}[ \E_{p,q,\sigma,\tau} [g(k_{t+1},l_{t+1},i_{t+1},j_{t+1}) | h_t ] ] \\
&=\sum_{t=0}^{T-1} \E_{p,q,\sigma,\overline{\tau}}[ \E_{p,q,\sigma,\overline{\tau}} [g(k_{t+1},l_{t+1},i_{t+1},j_{t+1}) | h_t ] ] =\gamma_T^{p,q}(\sigma,\overline{\tau})
\end{align*}
Moreover, since we proved that the law of $(h_t,l_{t+1})$ is the same under $\P_{p,q,\sigma,\tau}$ and under $\P_{p,q,\sigma , \overline{\tau}}$, this implies that $q_t(q,\tau)$ and $q_t(q,\overline{\tau})$ also have the same law, which proves the last assertion. \hfill $\Box$

\begin{rem} \rm
Note that in the above construction, the laws $\P_{p,q,\sigma,\tau}$ and  $\P_{p,q,\sigma , \overline{\tau}}$ do not coincide on the whole set $\Omega$ but generate the same payoffs. The reader may convince himself by considering the following example:
\[ L=\{1,2\}, \; N= \begin{pmatrix} \frac{1}{3} & \frac{2}{3} \\ \frac{2}{3} & \frac{1}{3} \end{pmatrix}, \; q=(\frac{1}{2},\frac{1}{2}), \; J=\{1,2\}, \]
and $\tau$ is defined by: play $1$ at stage $1$ and then at stage $t\geq 2$, play $1$ if $l_t=l_{t-1}$ and $2$ otherwise. A direct computation leads to 
\[ \P_{p,q,\sigma,\tau}( l_1=l_2=1, j_2=1)=\frac{1}{6} ,\quad \P_{p,q,\sigma , \overline{\tau}}( l_1=l_2=1, j_2=1)=\frac{1}{18} .\]
\end{rem}

\vspace{0,2cm}

The next proposition is a major point of our proof.\\

\begin{pro}\label{vhatequicontinu} 
For all $T \geq 1$, $\hat{v}_T$ is $3$-Lipschitz on $\Delta(K)\times \Delta(L)$. 
\end{pro}

\vspace{0,2cm}

The proof of Proposition \ref{vhatequicontinu}  is postponed to the Appendix. The crucial consequence  is to obtain uniform equi-continuity  of the family of functions $(\hat{v}_T)_{T\geq 1}$. In Renault \cite{R2006}, the functions $(\hat{v}_T)_{T\geq 1}$ were only proved to be u.s.c. and this property would not have been sufficient here, since the proof of our main Theorem \ref{thm1} relies on uniform convergence of the $(\hat{v}_T)_{T\geq 1}$ to their limit. The main point is to introduce a function $S(p,p')$ (see Definition \ref{definitionS}) which will play the role of a metric and with respect to which all the functions $(\hat{v}_T)_{T\geq 1}$ are ``non expansive''. \\

Before proving that the non revealing values 
$\hat{v}_T$ converge when $T$ goes to infinity,  we now present a      definition and a lemma,  allowing to concatenate non revealing strategies defined on blocks.
\begin{defi} \label{concat}
Given $T_0,T_1 \geq 1$ and a strategy $\sigma \in \Sigma$ such that for all $1 \leq t \leq T_1$, $\sigma_{T_0+t}$   does not depend on $(k_1,...,k_{T_0})$, we  define for all $h_{T_0} \in H_{T_0}$ a $T_1$-stage  strategy $\sigma(h_{T_0})$. This strategy  $\sigma(h_{T_0})$ plays,   after a history $(h'_{t-1},k'_1,...,k'_t)$,  what $\sigma$ plays at stage $T_0+t$ after  $((h_0, h'_{t-1}),k'_1,...,k'_t)$.
Formally,  
\[ \sigma(h_{T_0})_t (h'_{t-1},k'_1,...,k'_t) = \sigma_{T_0+t}((h_{T_0},h'_{t-1}),k'_1,...,k'_t), \forall 1 \leq t \leq T_1,   \forall h'_t \in H_{t-1}, \forall (k'_1,..,k'_t) \in K^t.\]
\end{defi}
 
\begin{lem}\label{concatenation}
Consider  $T_0,T_1 \in \N^*$ and strategies $(\sigma,\tau) \in \Sigma\times \mathcal{T}$ such that for all $1 \leq t \leq T_1$, $\sigma_{T_0+t}$ and $\tau_{T_0+t}$ do not depend on, respectively,  $(k_1,...,k_{T_0})$ and  $(l_1,..,l_{T_0})$.  The conditional law of $(k_{T_0+t},l_{T_0+t},i_{T_0+t},j_{T_0+t})_{t=1,..,T_1}$ given $h_{T_0}$ under $\P_{p,q,\sigma,\tau}$ is precisely $\P_{p_{T_0},q_{T_0},\sigma(h_{T_0}),\tau(h_{T_0})}$.
\end{lem}
The proof follows easily from a direct computation. 
\begin{rem}
Note that if $\tau \in \overline{\mathcal{T}}$, then the above condition on $\tau$ holds for all $T_0,T_1$.
We will use the following consequences in the sequel, at first
\[ \E_{p,q,\sigma,\tau} [ \sum_{t= 1}^{T_1} g(k_{T_0+t},l_{T_0+t},i_{T_0+t},j_{T_0+t}) ] =  \E_{p,q,\sigma,\tau} [ \gamma_{T_1}^{p_{T_0},q_{T_0}}(\sigma(h_{T_0}),\tau(h_{T_0})) ], \]
and for all  $t=1,..., T_1$ and any continuous function $f$ defined on $\Delta(K)\times \Delta(L)$, for all $h_{T_0}$ in $(I\times J)^{T_0}$,
$ \E_{p,q,\sigma,\tau}[ f(p_{T_0+t},q_{T_0+t}) |h_{T_0}]$$ =$$ \E_{p_{T_0},q_{T_0},\sigma(h_{T_0}),\tau(h_{T_0})}[f(p_t,q_t)],$
where $p_{T_0}=p_{T_0}(p, \sigma)(h_{T_0})$ and $q_{T_0}=q_{T_0}(q, \tau)(h_{T_0}).$
\end{rem}

 \vspace{0,5cm}
 
\begin{pro} \label{prohatv}$ \hat{v} =\lim_{T\to \infty} \hat{v}_T $ exists, is continuous on $\Delta(K)\times \Delta(L)$   and is   balanced. \end{pro}

 \vspace{0,5cm}
 
\noindent{Proof:}
Let $(p,q) \in \Delta(K)\times \Delta(L)$, $p^*=pB$, $q^*=qC$ and $\varepsilon>0$. Choose  $T_0,T_1 \in  \N^*$ such that 
$ \| M^{T_{0}} - B \| \leq \varepsilon,\quad \| N^{T_0}-C\| \leq \varepsilon$, $\quad \frac{T_0}{T_1} \leq \varepsilon, $
 and $\quad \hat{v}_{T_1}(p^*,q^*) \geq \lims_{T} \hat{v}_T(p^*,q^*) - \varepsilon. $ 
Given $N \in \N^*$ and an optimal strategy $\tau \in \overline{\mathcal{T}}$ of Player 2 in the game $\hat{\Gamma}_{N(T_0+T_1)} (p,q)$, 
let us construct the strategy $\sigma$ as follows. For $n=0$ to $N-1$, during stages $t=n(T_0+T_1) +1$ to $t=n(T_0+T_1) +T_0$, play a fixed action $i_0 \in I$. During the next $T_1$ stages, play an optimal strategy in the game $\hat{\Gamma}_{T_1}(p_{n(T_0+T_1) +T_0}(p,\sigma),q_{n(T_0+T_1) +T_0}(q,\tau))$. The payoff can be written as ($g_m$ being the payoff of stage $m$):

\begin{align*}
 N(T_0+T_1)\gamma^{p,q}_{N(T_0+T_1)}(\sigma,\tau) &= \sum_{n=0}^{N-1}\left( \sum_{m=n(T_0+T_1)+1}^{n(T_0+T_1)+T_0} \E_{p,q,\sigma,\tau}[ g_m]   + \sum_{m=n(T_0+T_1) + T_0 + 1}^{(n+1)(T_0+T_1)} \E_{p,q,\sigma,\tau}[g_m]  \right)\\
& \geq \sum_{n=0}^{N-1}\left( -T_0 + T_1\E_{p,q,\sigma,\tau}[ \hat{v}_{T_1}(p_{n(T_0+T_1) + T_0}(p,\sigma),q_{n(T_0+T_1) + T_0}(q,\tau))] \right)\\
& \geq \sum_{n=0}^{N-1}\left( -T_0 + T_1 \E_{p,q,\sigma,\tau}[ \hat{v}_{T_1}(p_{n(T_0+T_1)}(p,\sigma)M^{T_0},q_{n(T_0+T_1)}(q,\tau)N^{T_0})] \right)\\
& \geq \sum_{n=0}^{N-1}\left( -T_0 + T_1( \hat{v}_{T_1}(p^*,q^*) - 6\varepsilon ) \right).
\end{align*}

\noindent The first inequality follows from Lemma \ref{concatenation} and  the definition of $\sigma$. The second inequality is obtained by taking conditional expectation with respect to $h_{n(T_0+T_1)}$ and using that $\hat{v}_{T_1}$ is convex with respect to the second variable on $C^{-1}(q^*)$. The last inequality follows directly from the properties of $T_0$ and the fact that $\hat{v}_ {T_1}$ is $3$-Lipschitz. We deduce that for all $N \in \N^*$ 
\[ \hat{v}_{N(T_0+T_1)}(p,q) \geq \frac{ -T_0 + T_1( \hat{v}_{T_1}(p^*,q^*) - 6\varepsilon) }{(T_0+T_1)}\geq  \hat{v}_{T_1}(p^*,q^*) - 8\varepsilon \geq   \lims_{T} \hat{v}_T(p^*,q^*) - 9 \varepsilon.\]

\noindent It follows that 
$\limi_{N} \hat{v}_{N(T_0+T_1)}(p,q) \geq \lims_{T} \hat{v}_T(p^*,q^*) - 9 \varepsilon. $
Using that $\| \hat{v}_{T}- \hat{v}_{T+T'} \|_{\infty} \leq \frac{2T'}{T}$, this implies
  $\limi_{T} \hat{v}_{T}(p,q) \geq \lims_{T} \hat{v}_T(p^*,q^*)$.  Inverting the role of $K$ and $L$ leads to $\lims_{T} \hat{v}_{T}(p,q) \leq \limi_{T} \hat{v}_T(p^*,q^*)$,  
and we conclude that  $ \hat{v}(p,q)=\lim_{T\to \infty} \hat{v}_T(p,q)=\lim_{T} \hat{v}_T(p^*,q^*) $ exists and is a balanced function. It is continuous because  proposition \ref{vhatequicontinu} implies the uniform convergence of $ \hat{v}_T$. 
\hfill $\Box$

\section{The Mertens-Zamir system associated to the non-revealing value}\label{theoremMZ}
\p

We prove here an alternative formulation for $MZ(\hat{v})$ that we will use in the proof of theorem \ref{thm1} outlined in section 4.
Recall that we have assumed that the Markov chains on $K$ and $L$ are recurrent and aperiodic. 
We use here the notion of balanced function.
\begin{pro}\label{comparison}
$MZ(\hat{v})=\underline{v}= \overline{v}$ where 
 \[ \underline{v}= \inf  \left\{  w \in {\cal C}^+(\hat{v}) , w\;  balanced \right\}, \; \overline{v}=\sup   \left\{  w \in  {\cal C}^-(\hat{v}) , w \; balanced \right\} . \]
 \end{pro}

\vspace{0,3cm}

\noindent{Proof:} Consider  first any balanced function $w$ in $ {\cal C}$.  One may identify $w$ with its restriction $w^*$ to the product of simplices $P^* \times Q^*$. $I$-concavity and $\two$-convexity on $\Delta(K)\times \Delta(L)$ of $w$ is equivalent to   $I$-concavity and $\two$-convexity of $w^*$ on $P^*\times Q^*$. We will prove that  $\cav_I (w)$ is balanced and that $(\cav_I  w)^*= \cav_I(w^*)$.

Fix  $(p,q)$  in $\Delta(K)\times \Delta(L)$, and define $p^*=pB$ and $q^*=qB$. Recall that 
$\cav_{I} w(p,q) = \sup\{ \sum_{m=1}^N \alpha_m w(p_m,q) \;|\; N \geq 1, \forall m=1,..,N, \;\alpha_m \geq 0,\; \sum_{m=1}^N \alpha_m=1,\; \sum_{m=1}^N \alpha_m p_m =p \}.$
Given $(\alpha_m,p_m)_{m=1,..,N}$ as above, for all $m$, $w(p_m,q)=w(p_mB,q^*)$, and  $\sum_{m=1}^N \alpha_m p_mB =p^*$, which implies 
$ \sum_{m=1}^N \alpha_m w(p_m,q)  \leq \cav_{I} (w^*)(p^*,q^*), $
and therefore $\cav_{I} w(p,q) \leq \cav_{I} (w^*)(p^*,q^*)$. \\

  Define now the affine map $f$ from $P^*$ to $\Delta(K)$  by: 
\[ \forall p'  \in P^*, \forall r=1,..,r_M, \forall k \in K(r), \; f(p')^k= \left\{ \begin{matrix} (\sum_{s \in K(r)} p'^s)& \frac{p^k}{\sum_{s \in K(r)} p^s} & \text{if} & \sum_{s \in K(r)} p^s> 0 \\ (\sum_{s \in K(r)} p'^s)& p^*(r)^k & \text{if} & \sum_{s \in K(r)} p^s= 0 \end{matrix} \right. . \]
We have $f(p^*)=p$. Moreover,  for any recurrence class $K(r)$, we have $\sum_{k \in K(r)} f(p')^k= \sum_{s \in K(r)} p'^s$. This implies that  for all $p' \in P^*$, $f(p')B=p'$. Consider now  $(\alpha_m, p^*_m)_{m=1,..,N}$ such that
$ \forall m=1,...,N, \;\alpha_m \geq 0, p^*_m \in P^*,\; \sum_{m=1}^N \alpha_m=1,\; \sum_{m=1}^N \alpha_m p^*_m =p^*.$
We have by construction $\sum_{m=1}^N \alpha_m f(p^*_m) =p$ and $w(p^*_m,q^*)=w(f(p^*_m),q)$ which implies 
\[ \sum_{m=1}^N \alpha_m w(p^*_m,q^*)= \sum_{m=1}^N \alpha_m w(f(p^*_m),q) \leq \cav_{I} w(p,q),\]
and therefore $\cav_{I} (w^*)(p^*,q^*) \leq \cav_{I} w(p,q)$. We conclude that:
$\forall (p,q)\in \Delta(K)\times \Delta(L),\;\cav_{I} w(p,q)= \cav_{I} (w^*)(p^*,q^*).$
Consequently, $\cav_I (w)$ is balanced and  $(\cav_I  w)^*= \cav_I(w^*)$.\\

The functions $\underline{v}$ and $\overline{v}$ are clearly balanced, and we will consider $\underline{v}^*$ and $\overline{v}^*$. Since $\hat{v}$ is also balanced, we deduce that:
\begin{align*} 
\underline{v}^* & = \sup  \left\{ \begin{matrix} w: P^*\times Q^* \longrightarrow [-1,1], w\;  {\rm  continuous} \;\text{ $\two$-convex and}  \\
    \forall (p^*,q^*) \in P^*\times Q^*,\;  w(p^*,q^*) \leq \cav_{I} \Min(w,\hat{v}^*)(p^*,q^*) \end{matrix} \right\},
\end{align*}
and a similar property for $\overline{v}^*$. The equality $\underline{v}^*=\overline{v}^*$ follows from Theorem 2.1 in \cite{mertenszamir} and implies $\overline{v}=\underline{v}$. Moreover, the function $\overline{v}^*=\underline{v}^*$ is a solution of the Mertens-Zamir system on $P^* \times Q^*$:
 \begin{equation*} 
\forall (p',q') \in P^* \times Q^*,   \left\{ \begin{matrix} w(p',q') &=& \vex_{\two} \Max(w,\hat{v}^*)(p',q') \\
w(p',q') &=& \cav_{I} \Min(w,\hat{v}^*)(p',q') \end{matrix} \right. ,
\end{equation*}
For $(p,q)$ in $\Delta(K)\times \Delta(L)$ and $v=\overline{v}=\underline{v}$ , we have: 
\begin{eqnarray*}
v(p,q)=v^*(p^*, q^*) & = & \cav_{I} \Min(v^*,\hat{v}^*)(p^*,q^*)\\
 & =   &\left(\cav_{I} \Min(v,\hat{v}) \right)^* (p^*,q^*)\\
 & =   & \cav_{I} \Min(v,\hat{v} )   (p,q)
\end{eqnarray*}
Finally, $v$  is solution of the Mertens-Zamir system associated to $\hat{v}$.   \hfill $\Box$

\vspace{0,3cm}

\begin{rem}\label{cavtransient}
The above proof relies heavily on the fact that $B$ is associated to a recurrent Markov chain, that is no state is transient. Consider  for example 
the case where: $K=\{a,b,c\}$, and 
$ M= \begin{pmatrix} 1 &0 & 0 \\0 & 1 & 0 \\ \frac{1}{2} & \frac{1}{4}  & \frac{1}{4} \end{pmatrix} . $ 
The matrix $B=\lim_t M^t$ would     be  here the following:
$B=\begin{pmatrix} 1 &0 & 0 \\0 & 1 & 0 \\ \frac{2}{3} & \frac{1}{3} & 0 \end{pmatrix}.$ 
In this case, the function $f(p)=(p^b-1/3)^{+}$ (where $^+$ denotes the positive part) is such that $\cav f (p_0)=0$ and $\cav (f^*) (p_0B)=2/9$ for $p_0=(0,0,1)$. This was already pointed out in Renault \cite{R2006}, where general Markov matrices were considered. Our reduction of the problem to the study of aperiodic recurrent Markov chains allows therefore for useful technical simplifications. For example, in case $L$ is reduced to a single point, our characterization becomes $v(p)= \cav(\hat{v})(p)$. Using aperiodic but possibly not recurrent Markov chains, the characterization given in \cite{R2006}  reduces to  $v(p)=\cav(\hat{v})(pB)$ where $\hat{v}$ and $B$ are defined similarly. One deduces from the above proof that these two expressions coincide for recurrent Markov chains.
\end{rem}

\vspace{0,3cm}

\begin{rem} \label{conditioncav}
Note that one may add the condition that $w$ is also $I$-concave in the supremum defining   $\underline{v}$, because if $w$ fulfills the required assumptions, then  $\cav_{I}w$ also and $\cav_{I}w \geq w$. A similar result holds for $\overline{v}$.   
\end{rem}


We are now in a position to prove  theorem \ref{thm1}, that is to show that    $\lim_T v_T$ exists and   is $MZ(\hat{v}) = \inf\{  w \in {\cal C}^+(\hat{v})\}    = \sup   \{  w \in  {\cal C}^-(\hat{v})\}$. The proof is in the Appendix. 

\section{Open Questions}

 Note that the Lipschitz  constant of $3$  for the non revealing values is essentially used to obtain the uniform convergence  of the non revealing value functions $\hat{v}_T$. We do not know how to prove this uniform convergence in a simpler way, neither if one can obtain a better Lipschitz constant for the non revealing values.
 
The Maxmin of the infinitely repeated game with incomplete information on both sides was proved in \cite{aumasch} to be equal to $\cav_I \vex_{\two} u$. In the present model, it may be asked if the Maxmin is equal to  $\cav_I \vex_{\two} \hat{v}$. One may hope to   prove    that Player 2 can defend this quantity by combining the methods developed in the present work and the proof in \cite{mertenszamir}.  It should be more complex to determine what can be guaranteed by Player 1. Even if one proves that Player 1 can guarantee the limit $\hat{v}(p,q)$ {\it  in the   non revealing  game} (and consequently that the non revealing game has a uniform value), this does not imply a priori that Player 1 can guarantee $\vex_{II}\hat{v}$ in the original game. An other idea would be to consider the ``semi-revealing" game where only Player 1 is restricted to play a non revealing strategy, then  to prove that this semi-revealing game has a limit value $v'(p,q)$ that can be guaranteed by Player 1,  but then there is still little hope to be able to link $\cav_{I} v'$ with   $\vex_{\two} \hat{v}$.

 The results of Aumann and Maschler,  and of Mertens and Zamir  were extended to the case of imperfect observation of the actions. In the same way, our model can be extended to the case where instead of observing the past actions played, each player observes a stochastic signal depending only on the past actions (state independent signalling).  The notion of non-revealing strategies still makes sense, so conceptually one may hope to use the same ideas to prove the convergence of the value functions in this case. Note that a difficulty  in   this model  is that even if we fix the strategy of  a player, say Player 1,  this player is not able to compute along the play the belief of Player 2 about the current state in $K$. And   when introducing signals,  it is often the case that  one has to construct strategies by blocks sufficiently long to allow for statistical  tests. How these  blocks would interfer with the blocks where non revealing  strategies are used as in point (5) of the proof of Theorem \ref{thm1} is unclear, and this seems a   highly technical problem. 

  In the Mertens Zamir setup, as well for splitting games, the existence of the limit value has been extended to the case of general evaluations of the stream of payoffs (\cite{CLSsiam}), including the discounted approach. It seems natural to hope for the same kind of generalization here, however it is not clear how our proof could be extended to do so, and in particular to prove the convergence of the discounted values to the same limit $v$.

An important extension of repeated games with incomplete information on both sides  is the case of general, possibly correlated,  probabilities for the initial pair of states $(k,l)$. The results of Mertens and Zamir \cite{mertenszamir} were actually stated within this more general framework. This suggests the following generalization of our model: we are given an initial probability $\pi \in \Delta(K\times L)$ and a Markov chain $z_t=(k_t,l_t)$ with values in $K \times L$. At each stage, Player 1 observes $k_t$ and Player $2$ observes $l_t$. It is still possible to define the appropriate notion of convex and concave hull as  done by Mertens and Zamir. However, the definition of a non-revealing strategy should  be adapted. Indeed, even if Player 1 does not use his information on the state, the beliefs of Player 2 on the recurrence classes of the chain $z_t$ may evolve. This is left for future research. 

Let us finally mention the  recent work \cite{CRRV}, where a model of Markov Games with incomplete information is introduced, with transitions depending on a stage length parameter. The asymptotic analysis is of different nature and leads to the existence of a continuous-time limit value as the players play more and more frequently. This value is different of the long-time limit value considered here and the study of possible relationships between these two approaches is an interesting direction of research which is beyond the scope of this paper. 
 
\section{Appendix}

\subsection{Proof of Proposition \ref{vhatequicontinu}}\label{proofvhat}

Let us begin with a few definitions.  Recall that $K$ is partitioned into recurrence classes $K(1)$, ..., $K(r_M)$, each recurrent class $r$   having a unique invariant measure $p^*(r)\in \Delta(K(r))$. $P^*=\conv\{p^*(1),...,p^*(M)\}$ is the set of invariant measures   for $M$, as well as for  $B=\lim_t M^t$. In the following we denote by  $\Delta(r_M)$   the set of probabilities over the recurrence classes of $M$, i.e. over the finite set $\{1,...,r_M\}$. 
\begin{nota} Let $p$ be in $ \Delta(K)$.  For each  $r$, we denote by $\lambda(p)^r=\sum_{k \in K(r)} p^k$ the probability of the recurrence class $r$ under $p$. And we define 
the probability   $\lambda(p) \in \Delta(r_M)$    by $\lambda(p)=(\lambda(p)^r)_{r =1,...,r_M}.$
And for each $r=1,..,r_M$, we denote the conditional law, under $p$,    of $k$ given $r$ as the vector $p_{|r} \in \Delta(K)$ such that: 
\[ p_{|r}^k = \left\{ \begin{matrix} \frac{p^k}{\lambda(p)^r} & \text{if} & \lambda(p)^r > 0\,,\, k \in K(r), \\ p^*(r)^k & \text{if} & \lambda(p)^r = 0,\,\, k \in K(r),  \\ 0 & \text{if} & k \notin K(r).\end{matrix} \right. \]   
It follows that for all $p\in \Delta(K)$, $p= \sum_{r=1}^{r_M} \lambda(p)^r p_{|r}$.  Since $p_{|r}^k=0$ for $k \notin K(r)$, $p_{|r}$ will be often assimilated to a vector in $\Delta(K(r))$.
\end{nota}
\vspace{0.3cm}
 \begin{defi}\label{definitionS}
The   map $S:\Delta(K)\times \Delta(K) \rightarrow \R_+$ is defined by:
\[ \forall p, p' \in \Delta(K),\;\;S(p,p') = \| \lambda(p) - \lambda(p') \|+  \sum_{r\,:\,\lambda(p)^r\lambda(p')^r>0} \; \lambda(p')^r \; \|p_{|r} -  p'_{|r} \|. \]
\end{defi}
The asymmetric map $S$ will play an important role in the sequel, and will be used as   a (quasi, semi-) metric. $S$ is clearly not a metric, but separates points (we have $S(p,p')=0$ if and only if $p=p'$) and enjoys the following useful properties.
\vspace{0.2cm}
\begin{lem}\label{propertiesS}
For all $p,p' \in \Delta(K)$,
\[  \|p-p'\| \leq S(p,p') \leq 3\|p-p'\|, \;S(pM,p'M) \leq S(p,p'), \] 
\[S(pB,p'B)=\|pB-p'B\|= \| \lambda(p)-\lambda(p')\|, \;  \text{and} \;\; S(p,p')=\|p-p'\| \; \text{if} \;\; pB=p'B.  \]
\end{lem}

Before   proving  the   lemma, we will present the   roadmap  of  the proof. \\

\begin{defi} ~\\
\noindent 1) For $\alpha \in [0,1]$, define the operator $\Phi_{\alpha}$ on  $\mathcal{C}$  by 
\[ \Phi_{\alpha}(f)(p,q)=\max_{x \in NR(p)} \min_{y \in NR(q)} \left( \alpha \; G(p,q,x,y) + (1-\alpha) \sum_{i\in I, j \in J} x(p)(i)y(q)(j) f( \vec{p}(x,i)M, \vec{q}(y,j)N)\right). \] 
 \noindent 2) Let  $\mathcal{F}$ be the subset of  functions $f$ in ${\cal C}$  satisfying: 
\begin{itemize}
\item  $ \forall p, p' \in  \Delta(K)$, $\forall q \in   \Delta(L)$, $ \; f(p,q)-f(p',q) \leq S(p,p')$,
\item  $ \forall p  \in  \Delta(K)$, $\forall q, q' \in   \Delta(L)$,  $\; f(p,q)-f(p,q') \leq S(q',q)$,
\item $\forall p^*,q^* \in P^*\times Q^*$, $f$ is $I$-concave and $\two$-convex on $B^{-1}(p^*)\times C^{-1}(q^*)$.
\end{itemize}
\end{defi}
\vspace{0,3cm}
 We   now present our strategy for the  proof of Proposition \ref{vhatequicontinu}. We have  $\hat{v}_{t+1}= \Phi_{\frac{1}{t+1}} (\hat{v}_t)$ for $t \geq 0$ (with the convention $\hat{v}_0=0$), and all the functions 
in $\mathcal{F}$ are $3$-Lipschitz using Lemma \ref{propertiesS}. As a consequence, to show    that $\hat{v}_t$ is 3-Lipschitz for each  $t$, it is  sufficient to show that for all $\alpha$, $\mathcal{F}$ is stable by $\Phi_{\alpha}$. 

To do so, the main  point is the following: given $p$ and $p'$  in $\Delta(K)$ and $x$ in $NR(p)$, one need to find $x'$ in $NR(p')$ such that $(p',x')$ is not too far  from  $(p,x)$. If $p$ and $p'$ belong to the same set $B^{-1}(p^*)$, i.e. if $p$ and $p'$ assigns the same weight to each recurrence class, this can be efficiently done using, within each simplex $\Delta(K(r))$,  a lemma by Laraki (\cite{Rida} and \cite{Rida2},  see here the Appendix at the end). If $pB\neq p'B$,    there exists a simple affine map to transfer probabilities in $B^{-1}(pB)$ to probabilities in $B^{-1}(p'B)$ (Lemma \ref{appliaffine}). In the general case, we need to combine both aspects, applying first the affine transformation then the Laraki splitting  (see Lemma \ref{setvaluedmapH}). In the computations, the expression $S(p,p')$ appears and allows to control the expressions.

After the proof of Lemma \ref{propertiesS}, all the rest  of this section  is devoted to the proof that: $ \forall \alpha \in [0,1], \; \; \Phi_{\alpha}(\mathcal{F})\subset \mathcal{F}.  $
 
\vspace{0,3cm}
 
\noindent{Proof of Lemma \ref{propertiesS}.}
\begin{align*}
\|p-p'\|
&= 
\sum_{r \,:\,\lambda(p)^r\lambda(p')^r>0} \, \sum_{k \in K(r)} |\lambda(p)^r p_{|r}^k - \lambda(p')^r p'^k_{|r} |  \\
 &\qquad+ \sum_{r\,:\,\lambda(p)^r =0\,,\,\lambda(p')^r>0}\,\lambda(p')^r +\sum_{r\,:\, \lambda(p)^r>0\,,\,\lambda(p')^r=0}\,  \lambda(p)^r\\
&\leq 
\sum_{r\,:\,\lambda(p)^r\lambda(p')^r>0}\,\sum_{k \in K(r)} \left( |\lambda(p)^r p_{|r}^k - \lambda(p')^r p_{|r}^k | +| \lambda(p')^r p_{|r}^k - \lambda(p')^r p'^k_{|r}|  \right) \\
&\qquad + \sum_{r \,:\,\lambda(p)^r =0\,,\,\lambda(p')^r>0} \lambda(p')^r +\sum_{r \,:\, \lambda(p)^r>0\,,\,\lambda(p')^r=0} \lambda(p)^r \\
&=\sum_{r=1}^{r_M} |\lambda(p)^r - \lambda(p')^r | +\sum_{r \,:\, \; \lambda(p)^r\lambda(p')^r>0} \lambda(p')^r\sum_{k \in K(r)} | p_{|r}^k - p'^k_{|r} | \qquad = \qquad S(p,p').
\end{align*}
On the other hand, 
\begin{align*}
S(p,p')   
&= \| \lambda(p)  - \lambda(p')  \|+  \sum_{r\,:\,\lambda(p)^r\lambda(p')^r>0}\lambda(p')^r \| p_{|r} -  p'_{|r} \| \\
&\leq \| \lambda(p)  - \lambda(p')  \|+ \sum_{r\,:\,\lambda(p)^r\lambda(p')^r>0}\,\sum_{k \in K(r)} \left( |\lambda(p')^r p^k_{|r} - p^k | +| p^k- p'^k |  \right) \\
&= \|\lambda(p) - \lambda(p')\| +\sum_{r \,:\, \; \lambda(p)^r\lambda(p')^r>0} \left(|\lambda(p')^r- \lambda(p)^r| + \sum_{k \in K(r)} | p^k-p'^k |\right)\\
&\leq 2\|\lambda(p) - \lambda(p')\| + \| p-p'\| \leq 3 \| p-p'\|.
\end{align*}
For the next inequality, note that $\lambda(pM)=\lambda(p)$ and $(pM)_{|r}=p_{|r}M$, so that using that $p \mapsto pM$ is non-expansive
\begin{align*}
S(pM,p'M) &= \| \lambda(p)-\lambda(p')\|+\sum_{r\,:\,\lambda(p)^r\lambda(p')^r>0} \lambda(p') \| (p_{|r}- p'_{|r})M \| \\
&\leq \| \lambda(p)-\lambda(p')\|+\sum_{r\,:\,\lambda(p)^r\lambda(p')^r>0} \lambda(p') \| (p_{|r}- p'_{|r}) \| = S(p,p'). 
\end{align*}
Finally, the  equalities are  easily  proved  by  direct computation. \hfill $\Box$

\vspace{0.3cm}

Fix $p$ in $\Delta(K)$ and $q$ in $\Delta(L)$,  and let $p^*=pB$ and $q^*=qC$. To any   $x \in \Delta(I)^K$, we     associate $z \in \Delta(K\times I)$ defined by $z(k,i)= p^k x^k(i)$, and    will use  a more abstract way to denote non revealing strategies and conditional probabilities. Let $\Delta_f(\Delta(K))$ denote the set of probabilities $\mu$ with finite support on $\Delta(K)$. Such a probability can be written as 
\[ \mu= \sum_{n=1}^N \alpha_n \delta_{p_n},\]
with $(\alpha_n,p_n)_{n=1,..,N} \in ([0,1]\times \Delta(K))^N$ such that $\sum_{n=1}^N \alpha_n = 1$ and $p= \sum_{n=1,..,N} \alpha_n p_n$. The mean  $m(\mu)$ is defined as $ \int_{\Delta(K)} p d\mu(p)$. 
This set is endowed with the usual weak$^*$  topology which is  in particular metrized by the  Wasserstein (Kantorovich-Rubinstein)  distance induced by the norm $\|.\|$ and denoted $d_W$. Recall the standard duality formula: $ \forall \mu,\nu \in \Delta_f(\Delta(K)),$
$$\; d_W(\mu,\nu) = \max\{ \int \phi d\mu - \int \phi d\nu \;| \; \phi \in Lip_1\} = \min\{ \int_{p,p' \in \Delta(K)} \|p-p'\| d\pi(p,p') \;|\; \pi \in P(\mu,\nu)\}, $$
where $Lip_1$ denotes the set of $1$-Lipschitz functions on $\Delta(K)$, and   $P(\mu,\nu)$   the set of probability distributions on $\Delta(K) \times \Delta(K)$ having $\mu,\nu$ for marginals.  $\Delta_f(\Delta(K))$ is endowed with the  convex order  $\leq$   defined by
$$ \mu \leq \nu \Longleftrightarrow \int \phi d\mu \leq \int \phi d\nu \; \text{for all l.s.c. convex\; }  \phi: \Delta(K) \rightarrow \R \cup\{+\infty\} .$$ 
For instance, given $p$ and $p'$ in $\Delta(K)$ the Dirac measure   $\delta_{(p+p')/2}$ is smaller than the average $1/2 \; \delta_p +1/2\;  \delta_{p'}$. $\mu \leq \nu$ implies $m(\mu)=m(\nu)$.  And for all $\mu,\nu$ with finite support $U$ and $V$, then $\mu \leq \nu$ if and only if there exists $F: U \rightarrow \Delta(V)$ such that:  $\forall u \in U,\;  \sum_{v \in V} v F(u)(v)=u$  and $ \sum_{u \in U} \mu(u)F(u) = \nu$.   This last   condition   (``martingale decomposition")  can   be seen as follows: $\mu$ is the law of some random variable $X_1$ with values in $\Delta(K)$, $\nu$ is the law of some random variable $X_2$ with values in $\Delta(K)$, and we have the martingale condition: $\E(X_2|X_1)=X_1$ (see \cite{blackwell1953}).

\begin{defi}
Let $\Psi : \Delta(K \times I) \rightarrow \Delta_f(\Delta(K))$ be the disintegration mapping defined by 
\[ \Psi(z) = \sum_{i \in I \,:\, z(i)> 0} z(i) \delta_{\vec{p}(z,i)} \]
with $ z(i)=\sum_{k \in K} z(k,i)$ and $\vec{p}(z,i)^k = \frac{z(k,i)}{z(i)}$ for $i$ such that $z(i)>0$.
\end{defi}
\noindent   Imagine    $(k,i)$ is selected according to $z$ and only $i$ is observed, then for each $i$ one can compute the conditional $\vec{p}(z,i)$ on $K$.  $\Psi(z)$ gives  the law of this posterior on $K$.  It is standard that 
$\Psi$ is continuous and convex 
(see e.g. \cite{R2012} Lemmas 4.16 and 4.17). Using the $z$ variable and the desintegration $\Psi$, the definition of $NR(p)$   now reads  as follows: 
\begin{eqnarray*}
  x \in NR(p) & \Longleftrightarrow & \forall k \in K, \sum_{i \in I }z(k,i)=p^k \; \text{and}\; \Psi(z) \in \Delta_f(B^{-1}(p^*)),\\
 & \Longleftrightarrow & \delta_p \leq \Psi(z) \; \text{and}\; \Psi(z) \in \Delta_f(B^{-1}(p^*)).
  \end{eqnarray*}

\begin{defi} We define  the set-valued maps:
\[ 
  \begin{array}{cccl} H: & \Delta(K)  & \rightrightarrows &  \Delta_f(\Delta(K))\\ &
p \in \Delta(K)  & \rightarrow  & H(p)= \{ \mu \in \Delta_f(B^{-1}(pB)) | \; m(\mu)=p \}, \\ 
&&&\\
 R:  & \Delta_f(\Delta(K))  & \rightrightarrows & \Delta(K \times I) \\&
\mu \in \Delta_f(\Delta(K)) &\rightarrow & R(\mu)=\{z \in \Delta(K\times I) \,:\, \Psi(z) \leq \mu\}.
\end{array}\]  \end{defi}
 \noindent We define similarly $H: \Delta(L)  \rightrightarrows   \Delta_f(\Delta(L))$   and  $R: \Delta_f(\Delta(L)   \rightrightarrows \Delta(L \times J)$ (keeping the same letters $H$ and $R$ for simplicity). $H(p)$ is the set of probabilities $\mu$ over $B^{-1}(pB)$ with mean $p$, equivalently it is the set of probabilities $\mu$ over $B^{-1}(pB)$ such that $\delta_p \leq \mu$.  Hence $H(p)$  is the set of  splittings at $p$ such that all posterior keep the same weight on each recurrence class.    The set $R(\mu)$ can be seen the set of ``strategies" $z$ with disintegration no more  informative than $\mu$. It is easy to show that:

$$\forall \mu \in H(p), \forall z \in R(\mu),\;  \Psi(z) \in H(p).$$

For all  $\alpha\in [0,1]$ and $f$,    $\Phi_{\alpha}(f)$ can now be written as follows:
\begin{align*} 
\Phi_{\alpha}(f)(p,q) 
&= &\sup_{z \;s.t.  \Psi(z) \in H(p)} \;    \inf_{w \;s.t.  \Psi(w) \in H(q)} \left( \alpha G(z,w) + (1-\alpha) \int_{\Delta(K)\times \Delta(L)} f( \tilde{p}M, \tilde{q}N)d\Psi(z)\otimes\Psi(w)(\tilde{p},\tilde{q}) \right)\\
&=& \sup_{\mu \in H(p)}\; \sup_{z \in R(\mu)}\;  \inf_{\nu \in H(q)}\;  \inf_{w \in R(\nu)} \left( \alpha G(z,w) + (1-\alpha) \int_{\Delta(K)\times \Delta(L)} f( \tilde{p}M, \tilde{q}N)d\Psi(z)\otimes\Psi(w)(\tilde{p},\tilde{q}) \right),
\end{align*}
where $G(z,w)= \sum_{(k,l,i,j)\in K\times L \times I \times J}  z(k,i)w(l,j) g(k,l,i,j)$.\\

We now  simplify the above expression.
\begin{pro} \label{pro17}
Define $V(\mu,\nu) =\sup_{z \in R(\mu)}\inf_{w \in R(\nu)} G(z,w)$. Then $V$  is $1$-Lipschitz for $d_W$ and is $I$-concave and $\two$-convex on $\Delta_f(\Delta(K)) \times \Delta_f(\Delta(L))$. Moreover:
\[ \Phi_{\alpha}(f)(p,q)=\sup_{\mu \in H(p)} \inf_{\nu \in H(q)}  \left( \alpha V(\mu,\nu) + (1-\alpha) \int_{\Delta(K)\times \Delta(L)} f( \tilde{p}M, \tilde{q}N)d\mu\otimes\nu(\tilde{p},\tilde{q})\right). \]
\end{pro} 
 
\noindent{Proof of Proposition \ref{pro17}:}
Using that $f$ is $\two$-convex on $B^{-1}(p^*) \times C^{-1}(q^*)$ and that $C^{-1}(q^*)$ is invariant by $N$, we have for all $w \in R(\nu)$   
\[ \int_{\Delta(K)\times \Delta(L)} f( \tilde{p}M, \tilde{q}N)d\Psi(z)\otimes\Psi(w)(\tilde{p},\tilde{q}) \leq \int_{\Delta(K)\times \Delta(L)} f( \tilde{p}M, \tilde{q}N)d\Psi(z)\otimes\nu(\tilde{p},\tilde{q}), \]
which implies  the following equality
\begin{eqnarray*}
\; & \inf_{\nu \in H(q)} \inf_{w \in R(\nu)} \left( \alpha G(z,w) + (1-\alpha) \int_{\Delta(K)\times \Delta(L)} f( \tilde{p}M, \tilde{q}N)d\Psi(z)\otimes\Psi(w)(\tilde{p},\tilde{q})\right) \\
=& \inf_{\nu \in H(q)} \inf_{w \in R(\nu)} \left( \alpha G(z,w) + (1-\alpha) \int_{\Delta(K)\times \Delta(L)} f( \tilde{p}M, \tilde{q}N)d\Psi(z)\otimes\nu(\tilde{p},\tilde{q})\right)
\end{eqnarray*}
since one can always choose $\nu = \Psi(w)$ in the right-hand side.
A similar equality holds for $\Psi(z)$ and $\mu$. This implies that for the moment we have proved that 
for all $\alpha\in [0,1],$ $f$ continuous and $(p,q)\in \Delta(K)$,  \[ \Phi_{\alpha}(f)(p,q)=\sup_{\mu \in H(p)}\sup_{z \in R(\mu)} \inf_{\nu \in H(q)} \inf_{w \in R(\nu)} \left( \alpha G(z,w) + (1-\alpha) \int_{\Delta(K)\times \Delta(L)} f( \tilde{p}M, \tilde{q}N)d\mu\otimes\nu(\tilde{p},\tilde{q})\right). \]

Let us now study the properties of the set-valued map $R$.

\begin{cla}
The set-valued map $R: \mu \rightarrow R(\mu)$ is non expansive  from $(\Delta_f(\Delta(K)),d_W)$ to $(\Delta(K\times I), \|.\|)$.  Moreover,  for all $\beta \in [0,1]$, and $\mu,\mu' \in \Delta_f(\Delta(K))$, 
we have $ \beta R(\mu) + (1- \beta) R(\mu') \subset R(\beta \mu + (1-\beta) \mu'). $
Similar properties hold for $R: \nu \rightarrow R(\nu)$\ from $ \Delta_f(\Delta(L))$ to  $\Delta(L\times J)$. \end{cla}
 
\indent{Proof of the claim:}
We first prove the non expansive property. Fix $\mu$ and $\mu' \in \Delta_f(\Delta(K))$, and  $z \in R(\mu)$,  we have to find $z'$ in $R(\mu')$ such that $\|z-z'\|\leq d_W(\mu,\mu')$.  The idea is to  link   $\Psi(z)$ to $\mu$ via  the martingale decomposition, then to link $\mu$ and $\mu'$ via Kantorovich duality formula, and finally to define $z'$ using both links. 

Since $\Psi(z) \leq \mu$, there exists a map $F : \Delta(K) \rightarrow \Delta_f(\Delta(K))$ such that for all $x \in \Delta(K)$, $F(x)$ is centered in $x$ and $\mu = \sum_{i \in I} z(i)F(\vec{p}(z,i))$. 
Let $U,V$ denote respectively the supports of $\mu $ and $\mu'$, and consider by Kantorovich duality  a probability $\pi$ on $U \times V$ with     marginals $(\mu,\mu')$ such that $\sum_{(u,v)\in U\times V} \| u-v\| \pi(u,v) = d_W(\mu,\mu')$. For $u  \in U$, let $\vec{\pi}(u)  = (\frac{\pi(u,v)}{\mu(u)})_{v \in V}$ denote the conditional law on $V$ given  $u$ induced by $\pi$. Because  the second marginal  of $\pi$ is $\mu'$, we have $\sum_{u \in U} \mu(u) \vec{\pi}(u)= \mu'$. We define, for all $i$ in $I$:  
\[   \vec{p'}(i) = \sum_{v \in V}   v\;  m'(i)(v) \in \Delta(K)\;\text{with} \; m'(i)= \sum_{u\in U} F(\vec{p}(z,i))(u) \; \vec{\pi}(u) \in \Delta(V),  \] 
and $z'\in \Delta(K \times I)$ by 
$z'(k,i)=z(i)\; \vec{p'}(i)^k\;\;\;   \forall k\in K, \forall i\in I. $\\

We have $ \sum_{i \in I}z(i) m'(i)= \mu'$, and for each $i$ the probability $m'(i)$  is centered in $\vec{p'}(i)$. 
By construction 
$\Psi(z') = \sum_{i \in I} z(i) \delta_{\vec{p'}(i)}$. Using the martingale decomposition,  this implies    that  $\Psi(z')\leq \mu'$.  
Using Jensen's inequality, we have 
\begin{align*}
\| z - z'\| & = \sum_{i \in I} z(i)\| \vec{p}(z,i) - \vec{p'}(i) \|  
\\ &=   \sum_{i \in I} z(i)    \| \sum_{u \in U} uF(\vec{p}(z,i))(u) - \sum_{(u,v)\in U\times V} v F(\vec{p}(z,i))(u) \vec{\pi}(u)(v)  \| 
\\ &\leq \sum_{i \in I} z(i)  \sum_{u\in U}  F(\vec{p}(z,i))(u) \| u - \sum_{v\in V} v \vec{\pi}(u)(v)  \|
\\ &\leq \sum_{i \in I} z(i)  \sum_{u \in U} F(\vec{p}(z,i))(u)\sum_{v \in V} \vec{\pi}(u)(v) \| u-v \| 
\\ &= \sum_{u\in U} \sum_{v \in V} (\sum_{i \in I} z(i)   F(\vec{p}(z,i))(u)) \vec{\pi}(u)(v) \| u-v \| 
\\ &= \sum_{u\in U} \sum_{v \in V} \mu(u) \vec{\pi}(u)(v) \| u-v \|  = d_W (\mu,\mu')
\end{align*} 
The second assertion is due to the convexity of $\Psi$.
This ends the proof of the claim, and we now conclude the proof of Proposition \ref{pro17}.\\

 The law $\mu$ being fixed, let us  consider the map
\[ (z,\nu) \in R(\mu)\times H(q)  \rightarrow \inf_{w \in R(\nu)} \left( \alpha G(z,w) + (1-\alpha) \int_{\Delta(K)\times \Delta(L)} f( \tilde{p}M, \tilde{q}N)d\mu\otimes\nu(\tilde{p},\tilde{q})\right). \]
As an infimum of linear continuous functions, it is concave and upper semi-continuous with respect to $z$ on the compact convex set $R(\mu)$. 
On the other hand, the map 
\[\nu \rightarrow \int_{\Delta(K)\times \Delta(L)} f( \tilde{p}M, \tilde{q}N)d\mu\otimes\nu(\tilde{p},\tilde{q})   \]
is linear and using the second point in the above claim, the map 
$\left(\nu \rightarrow \inf_{w \in R(\nu)} G(z,w) \right)$
is convex  with respect to $\nu$.  It allows to apply Sion's Minmax Theorem (see \cite{sion}) to conclude that 
\[ \Phi_{\alpha}(f)(p,q)=\sup_{\mu \in H(p)} \inf_{\nu \in H(q)} \sup_{z \in R(\mu)}\inf_{w \in R(\nu)} \left( \alpha G(z,w) + (1-\alpha) \int_{\Delta(K)\times \Delta(L)} f( \tilde{p}M, \tilde{q}N)d\mu\otimes\nu(\tilde{p},\tilde{q})\right). \]
Using the claim about $R$ and that $z \rightarrow G(z,w)$ is $1$-Lipschitz, the function $V$ is $1$-Lipschitz for $d_W$. Moreover, using the same proof as above, $V$ is $I$-concave and $\two$-convex on $\Delta_f(\Delta(K)) \times \Delta_f(\Delta(L))$. 
The preceding expression becomes
\[ \Phi_{\alpha}(f)(p,q)=\sup_{\mu \in H(p)} \inf_{\nu \in H(q)}  \left( \alpha V(\mu,\nu) + (1-\alpha) \int_{\Delta(K)\times \Delta(L)} f( \tilde{p}M, \tilde{q}N)d\mu\otimes\nu(\tilde{p},\tilde{q})\right), \]
and this concludes the proof of Proposition \ref{pro17}. \hfill $\Box$
\vspace{0,3cm}

It remains now to study the properties of the set-valued map $H$. Let us start with the following lemma which proves a kind of ``triangle equality'' for the map $S$.
\begin{lem}\label{appliaffine}
For any pair $(p^*,p'^*)$ in $P^*$, there exists an affine map $L=L_{p^*,p'^*}$ from $B^{-1}(p^*)$ to $B^{-1}(p'^*)$   such that $L(p^*)=p'^*$, and for all $(p,p') \in   B^{-1}(p^*)\times B^{-1}(p'^*)$
\[  S(p,L(p))=\| L(p)-p \| =\|p^*-p'^*\| {\; \rm and \;}S(p,L(p)) + S(L(p),p') = S(p,p'). \]
\end{lem}
Note that using matrices of example A section 4, the map $L$ is just a homothetic transformation as illustrated below. 
\begin{center}
\includegraphics[scale=0.5]{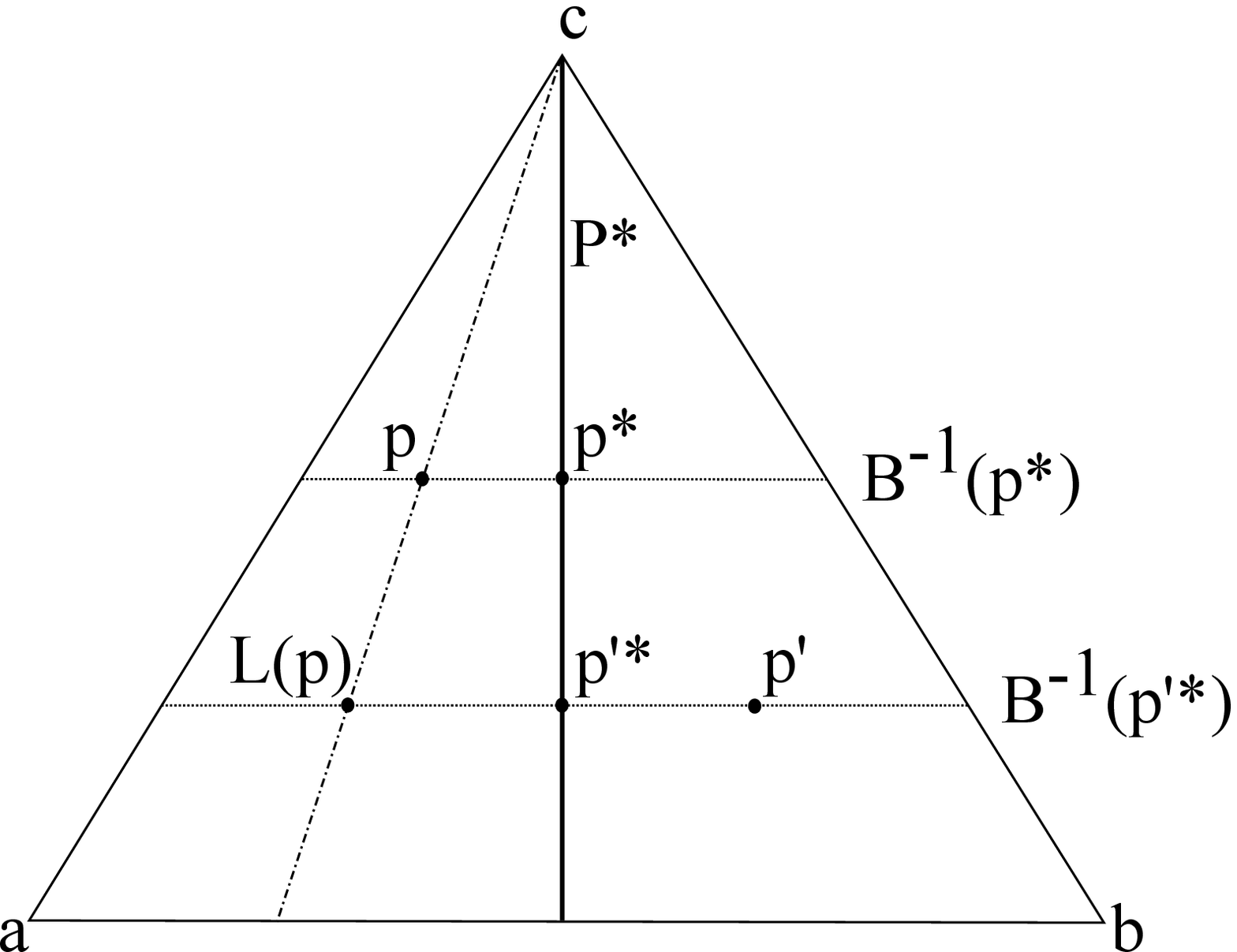}
\end{center}

\noindent{Proof:}
Note at first that $B^{-1}(p^*)= \{ \tilde{p} \in \Delta(K) | \lambda(\tilde{p}) = \lambda(p^*) \}$.  
Since $p^*,p'^* \in P^*$, we have   for all $r=1,..,r_M$, $p^*_{|r}=p'^*_{|r}=  p^*(r)$.  Define the map $L$ by 
\[  L(p) = \sum_{r=1}^{r_M} \lambda(p'^*)^r p_{|r}     \] 
It is easily seen that it defines an affine map with values in $B^{-1}(p'^*)$. Moreover,
\begin{align*}
 \| L(p)-p \| 
 &= \sum_{r\,:\, \lambda(p^*)^r>0} \sum_{k \in K(r)} | (\frac{\lambda(p'^*)^r}{\lambda(p^*)^r}-1) | p^k + \sum_{r \,:\,\lambda(p)^r=0} \lambda(p'^*)^r \sum_{k \in K(r)} p^*(r)^k  \\
&= \sum_{r \,:\, \lambda(p^*)^r>0} | \lambda(p'^*)^r-\lambda(p^*)^r| +\sum_{r \,:\, \lambda(p^*)^r=0}  \lambda(p'^*)^r \\
&= \sum_{r=1}^{r_M}  | \lambda(p'^*)^r-\lambda(p^*)^r| = \| p^*-p'^* \|= S(p^*,p'^*).
\end{align*}
For the last equality, note that by construction
$S(p,L(p))=  \|\lambda(p')-\lambda(p)\|$ and 
$S(L(p),p')=  \sum_{r\,:\,  \lambda(p')^r \lambda(p)^r > 0} \lambda(p')  \| p_{|r} - p'_{|r}\|,$ 
so that $S(p,L(p)) + S(L(p),p') = S(p,p')$. 
 \hfill $\Box$

We now recall a lemma   that will be used within each recurrence class.\footnote{The proof can be deduced from Lemma 8.2 in \cite{Rida},  which deals with general measurable spaces. A simpler  proof      can be  deduced from the version written in R. Laraki's PhD thesis \cite{Rida2}, proposition 5.12 page 107.} 
\begin{lem} \label{lemRida} (Laraki \cite{Rida}) ${\cal K}$ being a finite set, consider the  simplex $\Delta({\cal K})=\{(p^k)_{k\in {\cal K}}\in \R^{\cal K}_+, \sum_{k \in {\cal K}} p^k=1\}$ endowed with the $L^1$-norm $\|p-p'\|=\sum_{k \in {\cal K}} |p^k-p'^k|$. Consider a convex combination $p=\sum_{s \in S} \lambda_s p_s$ in $\Delta({\cal K})$, with $S$ finite, $\lambda_s\geq 0$ and $\sum_{s \in S} \lambda_s=1$. Then for every $p'$ in $\Delta({\cal K})$ there exists $(p'_s)_{s\in S}$ in $\Delta({\cal K})^S$ such that:$$p'=\sum_{s \in S} \lambda_s p'_s \;\; {\rm and}\;\; \sum_{s \in S} \lambda_s \|p_s-p'_s\|=\|p-p'\|.$$
\end{lem}

\noindent Using  the lemmas \ref{appliaffine} and \ref{lemRida}, we  now prove that   $H$ is ``non expansive  with respect to $S$".

\begin{lem}\label{setvaluedmapH}
 For all $p$ and $p'$ in $\Delta(K)$ and   for all $\mu=\sum_{n=1}^N \alpha_n \delta_{p_n} \in H(p)$,  there exists $\nu=\sum_{n=1}^N \alpha_n \delta_{p'_n}\in H(p')$ such that: 
\[  \sum_{n=1}^N \alpha_n S(p_n,p'_n) = S(p,p'). \]
In particular $p \rightarrow H(pB)$  is $1$-Lipschitz from $(\Delta(K), \|.\|)$ to   $(\Delta(\Delta(K)), d_{W})$.
\end{lem}

\noindent{Proof:}
Let $\mu= \sum_{n=1}^N \alpha_n \delta_{p_n} \in H(p)$ and $p' \in \Delta(K)$. Let us denote $p^*=pB$ and $p'^*=p'B$. 
Note that $p_{|r}= \sum_{n=1}^N \alpha_n p_{n|r}$ since for all $n=1,..,N$, $\lambda(p_n)=\lambda(p)$.
Let $L=L_{p^*,p'^*}$ be the affine map given by Lemma \ref{appliaffine}.
Define $\tilde{p}_n= L(p_n)$ for $n=1,..,N$, so that $L(p)= \sum_{n=1,..,N} \alpha_n \tilde{p}_n$. 
From Lemma \ref{appliaffine}, for all $p'' \in B^{-1}(p'^*)$, we have 
$S(p_n,p'') = S(p_n,\tilde{p}_n) + S(\tilde{p}_n,p''). $

For $r=1,..,r_M$, using Lemma \ref{lemRida} in the simplex $\Delta(K(r))$, there exists $(b_{r,n})_{n=1,..N} \in \Delta(K(r))$ such that 
$\sum_{n=1,..,N} \alpha_n b_{r,n}= p'_{|r}$ and 
$\sum_{n=1,..,N} \alpha_n \| \tilde{p}_{n|r} - b_{n,r}\| = \| \sum_{n=1,..,N} \alpha_n \tilde{p}_{n|r}-p'_{|r} \|.$
If $\lambda(p')^r> 0$, then 
$ \sum_{n=1,..,N} \alpha_n \tilde{p}_{n|r}=\sum_{n=1,..,N} \alpha_n p_{n|r}= p_{|r}. $

Let $p'_n = \sum_{r=1,..,r_M} \lambda(p')^r b_{r,n}$, then $\sum_{n=1}^N \alpha_n p'_n =p'$ and
\begin{align*} 
\sum_{n=1}^N \alpha_n S(\tilde{p}_n,p'_n) &= \sum_{n=1}^N \alpha_n \sum_{r \,:\, \lambda(p')^r > 0} \lambda(p')^r   \| \tilde{p}_{n|r}-p'_{n|r}\| 
\\ &= \sum_{n=1}^N \alpha_n \sum_{r \,:\, \lambda(p')^r > 0} \lambda(p')^r   \| \tilde{p}_{n|r}-b_{r,n}\|
\\&= \sum_{r \,:\, \lambda(p')^r > 0} \lambda(p')^r \| p_{|r}-p'_{|r} \|= S(L(p),p').
\end{align*} 
It follows that 
\[\sum_{n=1}^N \alpha_n S(p_n,p'_n)= \sum_{n=1}^N \alpha_n (S(p_n,\tilde{p}_n) +  S(\tilde{p}_n,p'_n)) 
= S(p,L(p)) + S(L(p),p') = S(p,p').\]
 \hfill $\Box$

\vspace{0,3cm}

We  now conclude the  general  proof of Proposition  \ref{vhatequicontinu}. 
Let $\varepsilon>0$, $p' \in \Delta(K)$ and  $\mu^* \in H(p)$ such that
\[\inf_{\nu \in H(q)}  \left( \alpha V(\mu^*,\nu) + (1-\alpha) \int_{\Delta(K)\times \Delta(L)} f( \tilde{p}M, \tilde{q}N)d\mu^*\otimes\nu(\tilde{p},\tilde{q}) \right)\geq \Phi_{\alpha}(f)(p,q)- \varepsilon.\]
If $\mu^*= \sum_{n=1}^N \alpha_n \delta_{p_n}$, using Lemma \ref{setvaluedmapH}, there exists $\mu'= \sum_{n=1}^N \alpha_n \delta_{p'_n} \in H(p')$ such that
$ \sum_{n=1}^N \alpha_n S(p_n,p'_n) = S(p,p'). $
Let $\nu^* \in H(q)$ such that 
\begin{multline*}
\alpha V(\mu',\nu^*) + (1-\alpha) \int_{\Delta(K)\times \Delta(L)} f( \tilde{p}M, \tilde{q}N) d\mu'\otimes\nu^*(\tilde{p},\tilde{q}) \\
\leq
\inf_{\nu \in H(q)}  \left( \alpha V(\mu',\nu) + (1-\alpha) \int_{\Delta(K)\times \Delta(L)} f( \tilde{p}M, \tilde{q}N)d\mu'\otimes\nu(\tilde{p},\tilde{q}) \right)+\varepsilon
\leq \Phi_\alpha(f) (p',q)+\varepsilon.
\end{multline*}
We have   for all $p,p'$, $f(pM,q)-f(p'M,q) \leq S(pM,p'M) \leq S(p,p')$, and  $d_W(\mu^*, \mu')\leq \sum_n \alpha_n \|p_n-p'_n\|$ $\leq \sum_N \alpha_n S(p_n,p'_n)$ $=S(p,p')$. Consequently, we have  
\begin{align*} 
\Phi_{\alpha}(f)(p,q) - \Phi_{\alpha}(f)(p',q)  
 & \leq  \alpha ( V(\mu^*,\nu^*)-V(\mu',\nu^*) ) +(1-\alpha)  \int f( \tilde{p}M, \tilde{q}N) d(\mu^*-\mu') \otimes\nu^*(\tilde{p},\tilde{q}) +2 \varepsilon 
\\ & \leq \alpha d_W(\mu^*,\mu') +(1-\alpha)  \int_{\Delta(L)} \sum_{n=1}^N \alpha_n(f(p_nM,\tilde{q}N) - f(p'_nM,\tilde{q}N))d\nu^*(\tilde{q})+2 \varepsilon 
\\ &\leq \alpha S(p,p') + (1-\alpha) \sum_{n=1}^N \alpha_n S(p_n,p'_n)+2 \varepsilon = S(p,p')+2 \varepsilon .
\end{align*}
We conclude that $\Phi_{\alpha}(f)(p,q) - \Phi_{\alpha}(f)(p',q)  \leq S(p,p')$. 
The symmetric property follows by the same method. In order to conclude the proof, it remains to prove that $\Phi_{\alpha}(f)$ is $I$-concave and $\two$-convex on $B^{-1}(p^*)\times C^{-1}(q^*)$, but this follows directly from the fact that for all $\beta \in [0,1]$ and $p,p' \in \Delta(K)$,  
$\beta H(p) + (1-\beta) H(p') \subset H(\beta p + (1-\beta)p'),$
and the similar property for $q,q' \in \Delta(L)$.

\subsection{Proof of Theorem \ref{thm1}} $\;$

The  proof of our  main theorem  \ref{thm1}, making formal the description  given in the end of section \ref{sec4} will use the two next lemmas. 
\begin{lem}\label{splitting}
Let $w \in {\cal C}$ be  such that 
$\forall (p,q) \in \Delta(K)\times \Delta(L), \; w(p,q) \leq \cav_{I} \Min(w,\hat{v})(p,q). $
For all $(p,q)$,  there exists a convex combination $(\alpha_m,p_m)_{m=1,..,K}$ in $ [0,1] \times \Delta(K)$ such that
$\sum_{m=1}^{K} \alpha_m = 1 , \; p=\sum_{m=1}^{K} \alpha_m p_m , \; \forall m=1,..,K, \; w(p_m,q) \leq \hat{v}(p_m,q),\; \text{and} \quad \sum_{m=1}^{K} \alpha_m w(p_m,q) \geq w(p,q).$
\end{lem}
\noindent The proof of lemma \ref{splitting} is standard: 
 there exists a convex combination $p=\sum_{m=1}^{K} \alpha_m p_m $ satisfying $\sum_{m=1}^{K} \alpha_m \Min(w,\hat{v})(p_m,q)= \cav_I \Min(w,\hat{v})(p,q).$ 
This implies that for all $m=1,..,K$, 
$ \Min(w,\hat{v})(p_m,q)= \cav_I \Min(w,\hat{v})(p_m,q) \geq w(p_m,q,$ 
so that $w(p_m,q)  \leq \hat{v}(p_m,q)$.
To conclude, note that
$\sum_{m=1}^{K} \alpha_m w (p_m,q) = \cav_I \Min(w,\hat{v})(p,q) \geq w(p,q). $ \hfill $\Box$ \\

\begin{lem} \label{SLemma2}
Consider  a convex combination $p=\sum_{s \in S}\alpha_s p_s$ in $\Delta(K)$.
Consider for each $s\in S$ a non revealing strategy $\sigma_s$ in $\hat{\Sigma}(p_s)$. There exists $\sigma$ in $\Sigma(p)$ such that: 
\[\forall q \in \Delta(L), \forall \tau \in {\cal T},\;\; \P_{p,q,\sigma,\tau}=\sum_{s \in S}  \alpha_s \P_{p_s,q,\sigma_s,\tau}.\]
As a consequence  we have for all $T \in \N^*$, $q$ in $\Delta(L)$ and $\tau\in {\cal T}$,
\[ \E_{p,q,\sigma,\tau}[ \sum_{t=0}^{T-1} \| \hat{q}_{t+1}(q,\tau)- \hat{q}_t(q,\tau) \| ] = \sum_{s \in S}  \alpha_s \E_{p_s,q,\sigma_s,\tau}[ \sum_{t=0}^{T-1} \| \hat{q}_{t+1}(q,\tau)- \hat{q}_t(q,\tau) \| ]. \]
Moreover, for any function $w \in \mathcal{C}$  which is $I$-concave, $\two$-convex  and all $t\geq 1$, we have 
$\E_{p,q,\sigma,\tau}[ w(\hat{p}_t,\hat{q}_t)] \geq \sum_{s\in S} \alpha_s w(p_sB,qC).$
\end{lem}

\vspace{0,3cm}

\noindent{Proof:}
Let us define $\sigma$ as in Lemma \ref{SLemma} via the splitting procedure: observe the first  state $k_1$ in $K$, then  choose a variable  $s$ in $S$ according  to the probability $\alpha_s p_s^{k_1}/p^{k_1}$, and play according  to $\sigma_s$. It defines a probability $\P$ on $\Omega \times S$ and using Kuhn's Theorem, there exists $\sigma$ such that the induced probability on $\Omega$ is equal to $\P_{p,q,\sigma,\tau}$. Moreover, for any event $A$ in $\Omega$, we have $\P(A|s)= \P_{p_s,q,\sigma_s,\tau}(A)$ almost surely, which implies the equality  $\P_{p,q,\sigma,\tau}=\sum_{s \in S}  \alpha_s \P_{p_s,q,\sigma_s,\tau}$. 

Recall that whenever $h_t$ is compatible with $(q,\tau)$, $q_t(q,\tau)(h_t)$ does not depend on $(p,\sigma)$. Fix $q$ and $\tau$, for simplicity  we write $q_t$ for $q_t(q,\tau)$, $\hat{q}_t$ for $\hat{q}_t(q,\tau)$ etc... Since $\hat{q}_t =q_t C$, the preceding equality implies that    
\[\E_{p,q,\sigma,\tau}[ \sum_{t=0}^{T-1} \| \hat{q}_{t+1} - \hat{q}_t  \| ] =\sum_{s \in S}  \alpha_s \E_{p_s,q,\sigma_s,\tau}[ \sum_{t=0}^{T-1} \| \hat{q}_{t+1} - \hat{q}_t  \| ].\]

Let us now work with the probability $\P$. Since $\sigma_s$ is non-revealing at $p_s$, we have for all stage $t$,   $s$ in $S$ and history $h_t$: 
$ (\P( k_{t+1} = k | h_t,s))_{k\in K} B =     p_t(p_s,\sigma_s)(h_t)B=p_s B. $
Denote by $\E$   the expectation under $\P$, and by $\tilde{s}$ the random variable with values in $S$ of the first choice of Player 1. The above implies: 
\[ \E[ p_{\tilde s}B | h_t ] = \E[ (\P( k_{t+1} = k | h_t,\tilde{s}))_{k\in K} | h_t]B = p_t(p, \sigma)(h_t)B =\hat{p}_t(p,\sigma)(h_t).\]
 Using that $w$ is $I$-concave and Jensen's inequality, we deduce that 
\[  \E[ w(p_{\tilde{s}}B, \hat{q}_t ) | h_t ] \leq \E[ w(\hat{p}_t(p,\sigma), \hat{q}_t ) | h_t ]. \]  
Using again that $\hat{q}_t$ does not depend on $(p,\sigma)$, we obtain 
\[ \E[ w(p_{\tilde{s}}B, \hat{q}_t  )]= \E [ \E[ w(p_{\tilde{s}}B, \hat{q}_t) | {\tilde{s}} ] ] = \sum_{s \in S} \alpha_s \E_{p_s,q,\sigma_s,\tau}[w(p_sB, \hat{q}_t  ]. \]
To conclude, $w$ being $\two$-convex, and by  Jensen's inequality again,  we get: 
\[\forall s \in S, \; \E_{p_s,q,\sigma_s,\tau}[w(p_sB, \hat{q}_t  ] \geq w(p_sB,qC). \]
\hfill $\Box$
\p
Let us turn to the main proof   of Theorem \ref{thm1}. It  is divided into  6 parts. 
\p
(1) {\it Fixing norms.} Recall that if $x=(x_s)_{s \in S}$ is an element of an euclidean space $\R^S$, $\|x\| = \sum_{s\in S} |x_s|$. For $T \in \N^*$, the set of strategies $\Sigma_T$ is the set of maps from $\cup_{t=0}^{T-1} (K \times I \times J)^{t-1} \times K$  to $\Delta(I)$ which is seen as a subset of $\R^I$. $\Sigma_T$ can be seen as a compact  subset of $(\R^I)^{\cup_{t=0}^{T-1} (K \times I \times J)^{t-1} \times K}$ endowed with the corresponding norm.
Note at first that for all $T \in \N^*$,  for all $(p,q) \in \Delta(K) \times \Delta(L)$, and for all $\sigma \in \Sigma_T$ and $\tau,\tau' \in \mathcal{T}_T$,
\[ | \gamma_T^{p,q}(\sigma,\tau)- \gamma_T^{p,q}(\sigma,\tau') |  \leq \| \tau-\tau' \|. \]

In the following points (2), (3), (4), the variables $T\geq 1$    and $\varepsilon\in (0,\frac{1}{I})$  are fixed. \\

\noindent (2) {\it Here we show that we can approximate strategies with low variations of the martingales $\hat{p}_t$ and $\hat{q}_t$ by non revealing strategies.}  Define $\Sigma_{T}^{\varepsilon}$ as the set of strategies $\sigma \in  \Sigma_{T}$ such that at any stage and for any history of the game, all the pure actions in $I$ are played with probability at least $\varepsilon$. 

Define the map $F$ from $\Delta(K)\times \Delta(L) \times \Sigma_{T} \times {\cal T}_{T}$ to $\R^{\vec{T}}$ by 
\[ F(p,q,\sigma,\tau) = \left(  \P_{p,q,\sigma,\tau} (h_{t+1}) (\hat{q}_{t+1}(q,\tau)(h_{t+1}) - \hat{q}_t(q,\tau)(h_{t}) )  \right)_{t=0,..,T-1\,,\, h_{t+1} \in (I\times J)^{t+1}}, \]
where for $h_{t+1}=(i_1,j_1,..,i_t,j_t, i_{t+1}, j_{t+1})$, $h_{t}$ is defined here as $(i_1,j_1,..,i_{t-1},j_{t-1}, i_t,j_t)$, and $\vec{T}$ is the large but finite integer $\vec{T}=L\sum_{t=0}^{T-1} (I \times J)^{t+1}$. Note that $F$ is continuous and  by construction 
$ \| F(p,q,\sigma,\tau ) \| = \sum_{t=0}^{T-1} \E_{p,q,\sigma,\tau}[ \|\hat{q}_{t+1}(q,\tau) - \hat{q}_{t}(q,\tau) \| ]. $

\indent Define the compact sets $R_{\varepsilon}=\Delta(K)\times \Delta(L) \times \Sigma_{T}^{\varepsilon} \times {\cal T}_{T}$
 and $\hat{R}_{\varepsilon}= F^{-1}(\{0\}) \cap R_\varepsilon$, $\hat{R}_{\varepsilon}$  is  the  subset of $R_{\varepsilon}$ consisting of elements  $(p,q,\sigma,\tau)$  such that $\tau$ is non-revealing at $q$.   
 
\begin{cla}  For all $\delta>0$, there exists $C_1=C_1(T,\varepsilon, \delta)$ such that 
\[ \forall (p,q,\sigma,\tau) \in R_{\varepsilon},\; \exists (\hat{p}, \hat{q},\hat{\sigma},\hat{\tau}) \in \hat{R}_{\varepsilon},\; \| (p,q,\sigma,\tau)-(\hat{p}, \hat{q},\hat{\sigma},\hat{\tau})\| \leq  \delta +C_1 \| F(p,q,\sigma,\tau) \|. \] \end{cla}
\noindent{Proof:} Assume by contradiction that for all $n \in \N^*$, there exists $(p_n,q_n,\sigma_n,\tau_n) \in R_{\varepsilon}$ such that for all $(\hat{p}, \hat{q},\hat{\sigma},\hat{\tau}) \in \hat{R}_{\varepsilon}$, 
$ \| (p_n,q_n,\sigma_n,\tau_n)-(\hat{p}, \hat{q},\hat{\sigma},\hat{\tau})\| >  \delta +n \| F(p_n,q_n,\sigma_n,\tau_n)\|.$
By compactness, we can extract a convergent subsequence with limit $(p^*,q^*,\sigma^*,\tau^*)$ such that $F(p^*,q^*,\sigma^*,\tau^*)=0$ and $ \| (p^*,q^*,\sigma^*,\tau^*)-(\hat{p}, \hat{q},\hat{\sigma},\hat{\tau})\| \geq  \delta.$
which is a contradiction. \hfill $\Box$\\
\p
\noindent (3) {\it Here we show that an optimal strategy of Player 1 in a non revealing game $\hat{\Gamma}_{T}(p,q)$ is good against strategies of Player 2 that are non revealing at   $q'$, where $q'$ is close to $q$}.

\begin{cla}  There exists $C_2=C_2(T,\varepsilon)$ such that for all $(p,q) \in \Delta(K)\times \Delta(L)$, for all $\sigma \in \hat{\Sigma}_{T}(p)$ which is optimal in $\hat{\Gamma}_{T}(p,q)$, for   all $q' \in \Delta(L)$, and for all $\tau' \in \hat{\mathcal{T}}(q')$, 
\[ \gamma_{T}^{p,q} (\sigma,\tau') \geq -\varepsilon - C_2 \|q-q'\| + \hat{v}_{T}(p,q) .\] \end{cla}

\vspace{0,3cm}
\noindent{Proof:} Assume by contradiction that for all $n\in \N^*$, there exists $(p_n,q_n) \in \Delta(K)\times \Delta(L)$, $\sigma_n \in \hat{\Sigma}_{T}(p_n)$ which is optimal in $\hat{\Gamma}_{T}(p_n,q_n)$,  $q'_n \in \Delta(L)$, and $\tau'_n \in \hat{\mathcal{T}}(q'_n)$ such that
$\gamma_{T}^{p_n,q_n} (\sigma_n,\tau'_n) \leq -\varepsilon - n \|q_n-q'_n\| + \hat{v}_{T}(p_n,q_n).$
By compactness, there exists a convergent subsequence with limit $(p^*,q^*,\sigma^*,q'^*,\tau'^*)$. Since the left hand-side of the above inequality is bounded below by $-1$, we have $q^*=q'^*$. Moreover, the  non-revealing graph  being closed, $\sigma^*$ is non-revealing at $p^*$ and $\tau'^*$ is non-revealing at $ q^*$. Using step (1), we have for all $n$,
\[ \hat{v}_{T}(p_n,q_n) -  \| \tau'^* - \tau'_n\| \leq \gamma_{T}^{p_n,q_n} (\sigma_n,\tau'^*)-  \| \tau'^* - \tau'_n\| \leq  \gamma_{T}^{p_n,q_n} (\sigma_n,\tau'_n). \]
It follows that  
$ \hat{v}_T(p^*,q^*) \leq   -\varepsilon + \hat{v}_{T}(p^*,q^*)$,
a contradiction. \hfill $\Box$\\

\noindent (4) {\it Here we control the error  while  perturbing an optimal strategy in  $\hat{\Gamma}(p,q)$ by a non revealing completely mixed strategy}.  

Let $\sigma^*_T$ denote the strategy in $\Sigma_T$ which plays an uniform distribution over $I$ at all stages independently of the history of the game, notice that $\sigma_T^*$ is non-revealing at all points in $\Delta(K)$.  
For all $(p,q) \in \Delta(K)\times \Delta(L)$, let $\sigma_T(p,q)$ denote an optimal strategy of Player 1 in $\hat{\Gamma}(p,q)$. Denote by   $\sigma_{T,\varepsilon}(p,q)$ the strategy which plays: $\sigma^*_T$ with probability $\varepsilon$, and $\sigma_T(p,q)$ with probability $1-\varepsilon$.  The strategy  $\sigma_{T,\varepsilon}(p,q)$ is non revealing at $p$, and we have  for all $\tau \in \mathcal{T}_T$,
$ \P_{p,q,\sigma_{T,\varepsilon},\tau}= \varepsilon \P_{p,q,\sigma_{T}^*,\tau} + (1-\varepsilon) \P_{p,q,\sigma_{T}(p,q),\tau}. $
We choose $\delta>0$ such that $\delta(1+C_2(T, \varepsilon))\leq \varepsilon$ and put $C_1=C_1(T, \varepsilon , \delta)$. Using step  (2),  for all $(p,q) \in \Delta(K) \times \Delta(L)$, and all   $\tau \in \mathcal{T}_T$,  
\[ \exists (\hat{p}, \hat{q},\hat{\sigma},\hat{\tau}) \in \hat{R}_{\varepsilon},\; \| (p,q,\sigma_{T,\varepsilon}(p,q),\tau)-(\hat{p}, \hat{q},\hat{\sigma},\hat{\tau})\| \leq  \delta +C_1  \| F(p,q,\sigma_{T,\varepsilon}(p,q),\tau) \|. \]
Using then (1) and (3), 
\begin{align*}
\gamma^{p,q}_T(\sigma_{T,\varepsilon}(p,q),\tau) &= \varepsilon \gamma^{p,q}_T(\sigma_{T}^*,\tau) + (1-\varepsilon) \gamma^{p,q}_T(\sigma_{T}(p,q),\tau) \\
&\geq -2 \varepsilon -  \| \tau - \hat{\tau} \| +\gamma^{p,q}_T(\sigma_{T}(p,q),\hat{\tau})\\
&\geq  -2 \varepsilon -  \| \tau - \hat{\tau} \| - \varepsilon - C_2  \| q-\hat{q} \| + \hat{v}_T(p,q) \\
 &\geq \hat{v}_T(p,q) - 3 \varepsilon - (1+C_2 ) (\delta + C_1 \| F(p,q,\sigma_{T,\varepsilon}(p,q),\tau\|) \\
 &\geq \hat{v}_T(p,q) - 4 \varepsilon - (1+C_2 )   C_1 \| F(p,q,\sigma_{T,\varepsilon}(p,q),\tau\|.
\end{align*}
\p
\noindent (5) {\it Here we   prove that:  $\limi_T v_T \geq \underline{v}$.} 

Let  $w \in {\cal C}$ be a balanced $I$-concave, $\two$-convex function such that 
$ \forall (p,q) \in \Delta(K)\times \Delta(L), \; w(p,q) \leq \cav_{I} \Min(w,\hat{v})(p,q). $
Fix  $(p,q)\in \Delta(K)\times \Delta(L)$, we have to  show that $\limi_{T} v_{T}(p,q) \geq w(p,q)$ (recall Remark \ref{conditioncav} and Proposition \ref{comparison} for the definition of $\underline{v}$).
\\

Fix  $\varepsilon \in (0, \frac{1}{I})$ and choose $T_0 \in \N^*$ such that $\| \hat{v}_{T_0} - \hat{v} \|_{\infty} \leq \varepsilon$.  For simplicity, we will write $p_t$ for $p_t(p,\sigma)$, $q_t$ for $q_t(q,\tau)$ and define $C_1$, $C_2$ for $C_1(T_0,\varepsilon, \delta)$, $C_2(T_0, \varepsilon)$ constructed  in the previous points (2), (3) and (4). Given $N \in \N^*$ and an optimal strategy $\tau \in \overline{\mathcal{T}}$ of Player 2 in $\Gamma_{NT_0}(p,q)$, we define a  strategy $\sigma$ as follows. 
\\

The set of stages is divided into consecutive blocks  of length $T_0$. For $n \geq 0$, at the beginning of block $n$, i.e. at the beginning of stage $nT_0+1$, Player 1 determines his strategy for the block according to the random variables $p_{nT_0}$ and $q_{nT_0}$. Two cases may occur. 

If $\hat{v}(p_{nT_0}, q_{nT_0}) \geq w(p_{nT_0 }, q_{nT_0})$, then $\sigma$ play  during the next $T_0$ stages the non revealing strategy $\sigma_{T_0,\varepsilon}(p_{nT_0}, q_{nT_0})$. 

If  $\hat{v}(p_{nT_0}, q_{nT_0}) < w(p_{nT_0}, q_{nT_0})$, using Lemma \ref{splitting}, there exists a convex combination $p_{nT_0}=\sum_{m=1}^{K} \alpha_{n,m} p_{n,m}$  which depends on $(p_{nT_0}, q_{nT_0})$, such that
$\forall m=1,..,K, \; w(p_{n,m},q_{nT_0}) \leq \hat{v}(p_{n,m},q_{nT_0})$ and 
$ \sum_{m=1}^{K} \alpha_{n,m} w(p_{n,m},q_{nT_0}) \geq w(p_{nT_0},q_{nT_0}).$
In this case Player 1 plays at block $n$ the strategy   $\overline{\sigma}_{T_0,\varepsilon}(p_{nT_0},q_{nT_0}) \in \Sigma_{T_0}$ given by Lemma \ref{SLemma2}  and such that for all $\tau' \in \mathcal{T}_{T_0}$,
$\P_{p_{nT_0},q_{nT_0}, \overline{\sigma}_{T_0,\varepsilon}(p_{nT_0},q_{nT_0}), \tau'}$$=$$ \sum_{m=1}^{K} \alpha_{n,m} \P_{p_{n,m}, q_{nT_0}, \sigma_{T,\varepsilon}(p_{n,m}, q_{nT_0}),\tau'}.$
This ends the definition of $\sigma$. 
\\

We have to show that the payoff $\gamma_{NT_0}(\sigma,\tau)$ is large. It can be written as 
\begin{align*} 
\gamma_{NT_0}(\sigma,\tau)&= \frac{1}{NT_0}\sum_{n=0}^{N-1} \E_{p,q,\sigma,\tau}[  \sum_{t=nT_0+1}^{(n+1)T_0}\E_{p,q,\sigma,\tau}[g(k_t,l_t,i_t,j_t)| h_{nT_0}] ] \\
&= \frac{1}{N}\sum_{n=0}^{N-1} \E_{p,q,\sigma,\tau} [\gamma_{T_0}^{p_{nT_0},q_{nT_0}} (\sigma(h_{nT_0}),\tau(h_{nT_0}))] \\
\end{align*}
where $\tau(h_{nT_0}) \in \mathcal{T}_{T_0}$ and $\sigma(h_{nT_0}) \in \Sigma_{T_0}$ are the continuation strategies defined in Def.  \ref{concat}. 
\begin{cla} For all $n=0,..,N-1$ and $h_{nT_0}$ occurring with positive probability,
\begin{align*}
\gamma_{T_0}^{p_{nT_0},q_{nT_0}} (\sigma(h_{nT_0}),\tau(h_{nT_0})) \geq 
w(p_{nT_0},q_{nT_0}) - 4 \varepsilon - (1+C_2) C_1 \| F(p_{nT_0},q_{nT_0},\sigma(h_{nT_0}),\tau(h_{nT_0}) \|. 
\end{align*} \end{cla}
\noindent{Proof:} For the case $\hat{v}(p_{nT_0}, q_{nT_0}) < w(p_{nT_0}, q_{nT_0})$, we have $\sigma(h_{nT_0})=\overline{\sigma}_{T_0,\varepsilon}(p_{nT_0},q_{nT_0})$ and
\begin{align*}
\gamma_{T_0}^{p_{nT_0},q_{nT_0}}& (\overline{\sigma}_{T_0,\varepsilon}(p_{nT_0},q_{nT_0}),\tau(h_{nT_0})) = \sum_{m=1}^{K} \alpha_{n,m} \gamma_{T_0}^{p_{nT_0},q_{nT_0}}(\sigma_{T,\varepsilon}(p_{n,m}, q_{nT_0}),\tau(h_{nT_0})) \\
&\geq 
\sum_{m=1}^{K} \alpha_{n,m}   \hat{v}_{T_0}(p_{n,m},q_{nT_0}) -4 \varepsilon  - (1+C_2) C_1\sum_{m=1}^{K} \alpha_{n,m}\| F(p_{n,m},q_{nT_0},\sigma_{T,\varepsilon}(p_{n,m}, q_{nT_0}),\tau(h_{nT_0})  \| 
\\ &\geq w((p_{nT_0},q_{nT_0}) - 4 \varepsilon  - (1+C_2) C_1\| F(p_{nT_0},q_{nT_0},\overline{\sigma}_{T,\varepsilon}(p_{nT_0}, q_{nT_0})),\tau(h_{nT_0}) \|,
\end{align*}
where the last inequality follows from the Lemmas \ref{splitting} and \ref{SLemma2}.  The case $\hat{v}(p_{nT_0}, q_{nT_0}) \geq w(p_{nT_0}, q_{nT_0})$ is simpler. \hfill $\Box$
 
\begin{cla} $\forall n=1,..,N-1,\;\; \E_{p,q,\sigma,\tau}[w(p_{nT_0},q_{nT_0})] \geq w(p,q).$\end{cla}
\noindent{Proof:}   By induction, it is sufficient to prove that for $n=0,..,N-2$,
\[  \E_{p,q,\sigma,\tau}[w(p_{(n+1)T_0},q_{(n+1)T_0}) | h_{nT_0}] \geq w(p_{nT_0},q_{nT_0}). \]
Let us at first consider the case $\hat{v}(p_{nT_0}, q_{nT_0}) \geq w(p_{nT_0}, q_{nT_0})$. 
Since $\sigma(h_{nT_0})$ is non-revealing at $p_{nT_0}$, we have  $p_{(n+1)T_0}=p_{nT_0}$. Then, using that $w$ is balanced, it is sufficient to prove that  
\[  \E_{p,q,\sigma,\tau}[w(p_{nT_0}B,q_{(n+1)T_0}C) | h_{nT_0}] \geq w(p_{nT_0}B,q_{nT_0}C). \]
The process $\hat{q}_{t}=q_tC$ being a martingale, the above inequality  follows from the fact that $w$ is $\two$-convex and Jensen's inequality.
Let us now consider the case $\hat{v}(p_{nT_0}, q_{nT_0}) < w(p_{nT_0}, q_{nT_0})$. Using Lemma \ref{SLemma2}, we have 
\begin{align*}
\E_{p,q,\sigma,\tau} &[w(p_{(n+1)T_0}B,q_{(n+1)T_0}C) | h_{nT_0}] \\
&= \E_{p_{nT_0},q_{nT_0},\overline{\sigma}_{T,\varepsilon}(p_{nT_0},q_{nT_0}),\tau(h_{nT_0})} [ w(\hat{p}_{T_0}(p_{nT_0},\overline{\sigma}_{T,\varepsilon}(p_{nT_0},q_{nT_0})),\hat{q}_{T_0}(q_{nT_0},\tau(h_{nT_0}))] \\
&\geq \sum_{m=1}^{K}\alpha_{n,m }w(p_{n,m}B,q_{nT_0}C) = \sum_{m=1}^{K}\alpha_{n,m }w(p_{n,m},q_{nT_0}) \geq w(p_{nT_0},q_{nT_0}). 
\end{align*}
 \hfill $\Box$
 
Summing up, we obtain:  
 \[\gamma_{NT_0}(\sigma,\tau)   \geq w(p,q) - 4 \varepsilon  
 - (1+C_2) C_1\frac{1}{N}\sum_{n=0}^{N-1} \E_{p,q,\sigma,\tau} [ \| F(p_{nT_0},q_{nT_0},\sigma(h_{nT_0}),\tau(h_{nT_0}))\|].\]
Finally, using Lemma \ref{concatenation}, and the classical bound on the $L_1$-variation of martingales (see for instance Proposition 3.8. in \cite{zamir}), we obtain:
\begin{align*} 
\sum_{n=0}^{N-1} \E_{p,q,\sigma,\tau} [ \| F(p_{nT_0},q_{nT_0},\sigma(h_{nT_0}),\tau(h_{nT_0})) \| ] 
&= \sum_{n=0}^{N-1}\E_{p,q,\sigma,\tau} [\E_{p,q,\sigma,\tau} [  \sum_{t=nT_0+1}^{(n+1)T_0} \| \hat{q}_{t} - \hat{q}_{t-1} \| | h_{nT_0}]]
\\&= \sum_ {t=1}^{NT_0} \E_{p,q,\sigma,\tau} [ \| \hat{q}_{t} - \hat{q}_{t-1} \| ]  \leq \sqrt{NT_0(L-1)}.
\end{align*}
We conclude that: $ 
\gamma_{NT_0}(\sigma,\tau)  \geq w(p,q) - 4 \varepsilon - (1+C_2) C_1\frac{\sqrt{T_0(L-1)}}{\sqrt{N}}.
$
Since $\tau$ is optimal in $\Gamma_{NT_0}(p,q)$, we deduce that 
$v_{NT_0}(p,q)  \geq w(p,q) -4 \varepsilon - (1+C_2) C_1\frac{\sqrt{T_0(L-1)}}{\sqrt{N}},$
and then 
$$ \limi_{N} v_{NT_0}(p,q) \geq w(p,q) - 4 \varepsilon . $$
Using that $\| v_{T}- v_{T+T'} \|_{\infty} \leq \frac{2T'}{T}$, this implies
$ \limi_{T} v_{T}(p,q) \geq w(p,q) - 4 \varepsilon$, 
and finally $\limi_{T} v_{T}(p,q) \geq w(p,q)$. We deduce that $\limi_T v_T(p,q) \geq \underline{v}$.\\

\noindent (6) {\it We conclude the proof of Theorem \ref{thm1}.} Symmetrically, we obtain as in point (5) that $\lims_T v_T(p,q) \leq \overline{v}$. The conclusion follows therefore from Proposition \ref{comparison}. \hfill$\Box$

\section*{Acknowledgments.}
  The authors gratefully acknowledge  the support of the Agence Nationale de la Recherche, under grant ANR JEUDY, ANR-10-BLAN 0112.

\end{document}